\documentclass[leqno,12pt]{article}
\usepackage{amsmath,amssymb,rotating,a4wide,color}
\usepackage[dvipsnames]{xcolor}
\usepackage{wasysym}
\usepackage{hyperref}
\usepackage{adjustbox}

\usepackage[all,cmtip]{xy}

\usepackage{verbatim}

\newdir^{ (}{{}*!/-3pt/\dir^{(}}
\newdir_{ (}{{}*!/-3pt/\dir_{(}}

\entrymodifiers={+!!<0pt,\fontdimen22\textfont2>}

\def\Agemotext{\rotatebox[origin=c]{180}{$\textstyle\kern-1pt\varOmega$}}
\def\Agemoscript{\rotatebox[origin=c]{180}{$\scriptstyle\kern-1pt\varOmega\kern1pt$}}
\def\Agemoscriptscript{\rotatebox[origin=c]{180}{$\scriptscriptstyle\kern-1pt\varOmega\kern1pt$}}
\def\Agemo{{\mathchoice{\Agemotext}{\Agemotext}{\Agemoscript}{\Agemoscriptscript}}}

\def\equaldown{\rotatebox[origin=c]{90}{$\textstyle\kern0pt=$}}

\def\strong{{\rm strong}}
\def\gr{\mathop{\rm gr}\nolimits}
\def\Sch{\mathop{\rm Sch}\nolimits}

\def\bigtimesdisplay{\mathop{\raise-2pt\hbox{\huge$\times$}}}
\def\bigtimestext{\mathop{\raise-1pt\hbox{\Large$\times$\kern-2pt}}}

\let\oldbigwedge\bigwedge
\def\newbigwedge{\mathord{\adjustbox{valign=B,totalheight=8.5pt}{$\oldbigwedge$}}}
\renewcommand{\bigwedge}{\newbigwedge}

\newbox\circbulletbox
\setbox\circbulletbox=\hbox{\,{\large$\circledcirc$}\kern-8.7pt\raise .6pt\hbox{$\bullet$}\,}

\let\le\leqslant
\let\ge\geqslant
\let\leq\leqslant

\let\triangleleftnaked\triangleleft
\def\triangleleft{\mathrel{\triangleleftnaked}}
\def\depsilon{d^{\kern1pt\epsilon}}
\def\Fratt{\mathop{\rm Fr}\nolimits}
\def\cond{{\rm cond}}
\def\obs{{\rm obs}}
\def\circVbig{\hbox{\text{\it\r{V}}}}
\def\circVscript{\hbox{\scriptsize\text{\it\r{V}}}}
\def\circVscriptscript{\mbox{\tiny\text{\it\r{V}}}}

\def\circVlimits_#1^#2{{\mathchoice%
   {\circVbig{}^{\kern2pt #2}_{\kern-2pt #1}}%
   {\circVbig{}^{\kern2pt #2}_{\kern-2pt #1}}%
   {\scriptstyle\circVscript{}^{\kern1.7pt #2}_{\kern-1pt #1}}%
   {\scriptscriptstyle\circVscriptscript{}^{\kern1.5pt #2}_{\kern-1pt #1}}%
   }}
\def\circVr_#1{\circVlimits_#1^r}
\def\circVs_#1{\circVlimits_#1^s}

\def\circWbig{\hbox{\text{\it\r{W}}}}
\def\circWscript{\hbox{\scriptsize\text{\it\r{W}}}}
\def\circWscriptscript{\mbox{\tiny\text{\it\r{W}}}}

\def\circWlimits_#1^#2{{\mathchoice%
   {\circWbig{}^{\kern2pt #2}_{\kern-2pt #1}}%
   {\circWbig{}^{\kern2pt #2}_{\kern-2pt #1}}%
   {\scriptstyle\circWscript{}^{\kern1.7pt #2}_{\kern-1pt #1}}%
   {\scriptscriptstyle\circWscriptscript{}^{\kern1.5pt #2}_{\kern-1pt #1}}%
   }}

\def\OM{\mathchoice
  {\rlap{\kern3.2pt$\overline{\phantom{L}}$}M}
  {\rlap{\kern3.2pt$\overline{\phantom{L}}$}M}
  {\rlap{\kern2.4pt$\scriptstyle\overline{\phantom{L}}$}M}
  {\rlap{\kern1.8pt$\scriptscriptstyle\overline{\phantom{L}}$}M}}

\def\mycirc{{\kern1pt\circ\kern2pt}}

\def\tra{\mathop{\rm tra}\nolimits}
\def\Cl{\mathop{\mathcal Cl}\nolimits}

\def\Rel{\mathop{\rm Rel}\nolimits}

\def\Aut{\mathop{\rm Aut}\nolimits}

\def\Gal{\mathop{\rm Gal}\nolimits}

\def\Spec{\mathop{\rm Spec}\nolimits}

\def\Ker{\mathop{\rm Ker}\nolimits}

\def\Stab{\mathop{\rm Stab}\nolimits}

\def\Ad{\mathop{\rm Ad}\nolimits}

\def\GL{\mathop{\rm GL}\nolimits}
\def\SL{\mathop{\rm SL}\nolimits}

\def\PGL{\mathop{\rm PGL}\nolimits}

\def\ab{{\rm ab}}

\let\phi\varphi

\def\theta{\vartheta}
\let\epsilon\varepsilon
\let\setminus\smallsetminus

\newtheorem{Thm}{Theorem}[section]
\newtheorem{Prop}[Thm]{Proposition}
\newtheorem{Lem}[Thm]{Lemma}

\newtheorem{Def}[Thm]{Definition}

\newtheorem{Rem}[Thm]{Remark}

\numberwithin{Thm}{section}
\def\UseTheoremCounterForNextEquation{\setcounter{equation}{\value{Thm}}\addtocounter{Thm}{1}}

\def\qed{{\hskip0pt\unskip\unskip\nobreak\hfil\penalty50
          \hskip1em\hbox{}\nobreak\hfil
           {$\square$}
          \parfillskip=0pt\finalhyphendemerits=0
          \par}\medskip}

\newenvironment{Proof}
               {\noindent{\bf Proof.}\ }
               {\qed}

               {\noindent{\bf Proof of #1.}\ }
               {\qed}

\newcommand{\BF}{{\mathbb{F}}}

\newcommand{\BQ}{{\mathbb{Q}}}

\newcommand{\BZ}{{\mathbb{Z}}}

\newcommand{\Fg}{{\mathfrak{g}}}

\newcommand{\Fl}{{\mathfrak{l}}}

\newcommand{\Fs}{{\mathfrak{s}}}

\newcommand{\CG}{{\cal G}}

\newcommand{\CO}{{\cal O}}

\usepackage{bbm}
\usepackage{tikz-cd}
\usepackage{mathrsfs}
\usepackage{siunitx}

\newcommand{\sO}{{\mathcal O}}
\DeclareMathOperator{\HOM}{\kern -2pt\mathscr{H}\text{\kern -2.5pt {\emph{om}}}}

\newbox\mybox
\def\arrover#1{\mathrel{
       \setbox\mybox=\hbox spread 1.4em
              {\hfil$\scriptstyle#1$\hfil}
       \vbox{\offinterlineskip\copy\mybox
             \hbox to\wd\mybox{\rightarrowfill}}}}

\def\larrover#1{\mathrel{
       \setbox\mybox=\hbox spread 1.4em
              {\hfil$\scriptstyle#1\vphantom{g}$\hfil}
       \vbox{\offinterlineskip\copy\mybox
             \hbox to\wd\mybox{\leftarrowfill}}}}

\def\ontoover#1{\mathrel{
       \setbox\mybox=\hbox spread 1.4em
              {\hfil$\scriptstyle#1\vphantom{g}$\hfil}
       \vbox{\offinterlineskip\copy\mybox
             \hbox to\wd\mybox{\rightarrowfill\hskip-2.8mm
                               $\rightarrow$}}}}
\def\leftontoover#1{\mathrel{
       \setbox\mybox=\hbox spread 1.4em
              {\hfil$\scriptstyle#1\vphantom{g}$\hfil}
       \vbox{\offinterlineskip\copy\mybox
             \hbox to\wd\mybox{$\leftarrow$\hskip-2.8mm
                               \leftarrowfill}}}}
\let\longto\longrightarrow

\let\onto\twoheadrightarrow
\def\longonto{\ontoover{\ }}

\def\isoto{\mathrel{
       \setbox\mybox=\hbox spread 0.9em
              {\hfil$\scriptstyle\sim$\hfil}
       \vbox{\offinterlineskip\copy\mybox
             \hbox to\wd\mybox{\rightarrowfill}}}}

\def\Bigskip{\bigskip\bigskip}

\begin{document}

\title{Schur $\sigma$-groups of type $(3,3)$ for $p=3$}

\author{\qquad
\begin{minipage}{.3\hsize}
Eric Ahlqvist\\[12pt]
\small Department of Mathematics \\
Stockholm University\\
106 91 Stockholm\\
Sweden \\
eric.ahlqvist@math.su.se\\[9pt]
\end{minipage}
\qquad
\begin{minipage}{.3\hsize}
Richard Pink\\[12pt]
\small Department of Mathematics \\
ETH Z\"urich\\
8092 Z\"urich\\
Switzerland \\
pink@math.ethz.ch\\[9pt]
\end{minipage}
\qquad
}

\date{\today}
%\date{May 14, 2025}

\maketitle

\Bigskip

\begin{abstract}
For any imaginary quadratic field $K$, the Galois group $G_K$ of its maximal unramified pro-$3$-extension is a Schur $\sigma$-group. If this has Zassenhaus type $(3,3)$, there are 13 possibilities for the isomorphism class of the finite quotient $G_K/D_4(G_K)$. 

We prove that for 10 of these 13 cases $G_K$ is either finite or isomorphic to an open subgroup of a form of $\PGL_2$ over~$\BQ_3$. Combined with the Fontaine-Mazur conjecture, or with earlier work on an analogue of the Cohen--Lenstra heuristic for Schur $\sigma$-groups, this lends credence to the ``if'' part of a conjecture of McLeman.

Using explicit computations of triple Massey products, we also test the heuristic for all imaginary quadratic fields~$K$ with $d(G_K)=2$ and discriminant $-10^8 < d_K < 0$ and find a reasonably good agreement.
\end{abstract}

{\renewcommand{\thefootnote}{}
%11R11: Quadratic extensions
%11R32: Galois theory
%11R34: Galois cohomology
%20D15: Nilpotent groups, $p$-groups
%20E18: Limits, profinite groups
%20F05: Generators, relations, and presentations
%
%11R29: Class numbers, class groups, discriminants
%11R37: Class field theory
%11R45: Density theorems
%20F14: Derived series, central series, and generalizations
}

\newpage
\renewcommand{\baselinestretch}{0.6}\normalsize
\tableofcontents
\renewcommand{\baselinestretch}{1.0}\normalsize

%%%%%%%%%%%%%%%%%%%%%%%%
% EDITORIAL COMMENTS:
%At the end, check that we always say `weak Schur $\sigma$-groups' where necessary.\\
%Check for systematic use of $\dim_{\BF_p}(...)$ versus $\dim(...)$.\\
%Always say ``quotient" instead of ``factor group".\\

%%%%%%%%%%%%%%%%%%%%%%%%%%%%%%%%%%%%%%%%%%%%%%%%
%%%%%%%%%%%%%%%%%%%%%%%%%%%%%%%%%%%%%%%%%%%%%%%%
%\newpage
\section{Introduction}
\label{Intro}

Fix an odd prime number~$p$. For any imaginary quadratic field~$K$, the compositum $\widetilde{K}$ of all fields in the Hilbert $p$-class field tower of~$K$ is the maximal unramified pro-$p$ extension of~$K$. The problem of determining whether this is finite over $K$ goes back to Hilbert and Furtwängler in the 1920s, and is closely connected to the classical question of whether $K$ possesses a finite extension whose ring of integers is a unique factorization domain; see \cite{RoquetteCFT}. Despite significant progress, no complete characterization is known for the pairs $(p,K)$ for which $\widetilde{K}/K$ is finite.

The problem is equivalent to whether the Galois group $G_K := \Gal(\widetilde{K}/K)$, called the \emph{$p$-tower group of~$K$}, is finite. This is a finitely generated pro-$p$-group, and since $\widetilde{K}$ is Galois over~$\BQ$, it possesses an automorphism of order at most $2$ induced by complex conjugation. Axiomatizing these and other known properties of $G_K$ leads to the notions of \emph{(weak or strong) Schur $\sigma$-groups}, which can be constructed and studied abstractly and seem interesting in their own right.

%%%%%%%%%%%%%%%%%%%%%%%%
\medskip
So let $G$ be a Schur $\sigma$-group with minimal number of generators $d(G)$. Then by Golod--Shafarevich \cite{GolodShafarevich1964} and Koch--Venkov \cite{KochVenkov1974}, it is known that $G$ is finite if $d(G) \le 1$, infinite if $d(G) \ge 3$, and in the case $d(G) = 2$ it is infinite unless its Zassenhaus type is $(3,3)$, $(3,5)$, or $(3,7)$ (see McLeman \cite[\S2]{McLeman2008}). The remaining cases are not yet completely settled, but since the work of Scholz and Taussky \cite{ScholzTaussky1934}, it has been known that there exist pairs $(p,K)$ for which $G_K$ is finite of type~$(3,3)$. 
Based on this, and on the observation that the Golod--Shafarevich kind of arguments to prove that $G$ is infinite only barely fail in the cases $(3,5)$ and $(3,7)$, McLeman \cite[Conj.\,2.9]{McLeman2008} conjectured that a $p$-tower group $G_K$ is finite if and only if $d(G_K)\le1$ or $G_K$ has Zassenhaus type~$(3,3)$. 

The ``only if'' part of this conjecture has recently been refuted by the first author and Carlson \cite[Thm.\,5.7]{AhlqvistCarlson2025}, who found pairs $(p,K)$ where $G_K$ is finite of Zassenhaus type~$(3,5)$. In particular there therefore exist finite Schur $\sigma$-groups of Zassenhaus type $(3,3)$ and~$(3,5)$, which is not yet known for type $(3,7)$. 

The ``if'' part of McLeman's conjecture, however, remains open. The paper \cite{Pink2025} by the second author provides some evidence for it in the case $p>3$. Among other things it proves that any infinite strong Schur $\sigma$-group of Zassenhaus type~$(3,3)$ is isomorphic to an open subgroup of a form of $\PGL_2$ over~$\BQ_p$. As the unramified Fontaine--Mazur conjecture \cite[Conj.\,5b]{FontaineMazur1995} excludes this possibility for $p$-tower groups, the ``if'' part of McLeman's conjecture becomes a direct consequence of the Fontaine--Mazur conjecture for $p>3$.

%%%%%%%%%%%%%%%%%%%%%%%%
\medskip
In the present paper we investigate weak Schur $\sigma$-groups of Zassenhaus type $(3,3)$ for $p=3$. This case differs from the case $p>3$ in the following way. 

\medskip
Let $F_n$ denote the free pro-$p$-group on $n$ generators, so that any weak Schur $\sigma$-group $G$ with $d(G)=n$ is isomorphic to a quotient of~$F_n$. 
%For any pro-$p$-group $G$ let $P_j(G)$ and $D_j(G)$ denote the $j$-th steps of the lower $p$-central series, respectively the Zassenhaus filtration, with the convention that $P_0(G)=G$ and ${D_1(G)=G}$. 
%Recall that the Zassenhaus filtration of a finitely generated pro-$p$-group is a descending sequence of open normal subgroups $G=D_1(G)\supset D_2(G)\supset\ldots$. 
For any pro-$p$-group $G$ let $D_j(G)$ denote the $j$-th step of the Zassenhaus filtration, beginning with ${D_1(G)=G}$. Then the subquotient $\gr_3^D(F_2)$ of the Zassenhaus filtration of $F_2$ is an $\BF_p$-vector space of dimension~$2$ if $p>3$, and of dimension~$4$ if $p=3$. It follows that among all weak Schur $\sigma$-groups $G$ of Zassenhaus type $(3,3)$, for the isomorphism class of the finite quotient $G/D_4(G)$ there is only one possibility if $p>3$, while there are $13$ different possibilities if $p=3$ (see Section~\ref{33ByD4}).

%%%%%%%%%%%%%%%%%%%%%%%%
\medskip
%{\bf Schur $\sigma$-groups and $p$-adic analytic groups:}
For 10 of those 13 cases, we manage to prove the same results as in \cite{Pink2025}, though the technique is much more involved.

\medskip
Recall that a pro-$p$-group $G$ is called powerful if its Frattini subgroup is generated by $p$-th powers. By a theorem of Lazard \cite{Lazard1965}, \cite[Thm.\,8.1]{DdSMS2003}, a finitely generated pro-$p$-group has the structure of a $p$-adic analytic group if and only if it possesses a powerful open subgroup. A crucial step in \cite{Pink2025} was to show that the subgroup $D_3(G)$ is powerful. In the present case $p=3$ we are able to decide precisely when $D_2(G)$ is powerful, and we also obtain partial results for $D_3(G)$ and $D_4(G)$.

For this and for possible later use in analyzing the remaining cases and in other situations, we develop a machinery to decide whether a more general open normal subgroup $E(G)$ of $G$ is powerful. This machinery applies to weak Schur $\sigma$-groups $G$ with arbitrary $d(G)$ and arbitrary odd~$p$. 
By Proposition~\ref{EGisPowerful} the condition depends only on the finite quotient $G/E_2(G)$ for some canonically defined smaller open normal subgroup $E_2(G)$. This makes it possible to test whether $E(G)$ is powerful using a computer algebra system.

To verify the criterion for all weak Schur $\sigma$-groups $G$ with a given isomorphism class of $G/D_4(G)$, we need to construct all possibilities for $G/E_2(G)$. There are too many cases to do this by hand, so we use the computer algebra system GAP \cite{GAP4} that specializes in group theory.
Its ANUPQ package \cite{ANUPQ}, which is a very efficient implementation of O'Brien's $p$-group generation algorithm \cite{OBrien1990}, is particularly convenient for the construction of $G/E_2(G)$ by recursion on~$E_2$. 

In this way we show that $D_2(G)$ is powerful in 10 of the above 13 cases; see the annotated worksheet \cite[{\tt 1-Checking-powerfulness}]{Ahlqvist-Pink-I}. In particular $G$ is then a $p$-adic analytic group. If moreover $G$ is a strong Schur $\sigma$-group, which specifically holds for any $p$-tower group~$G_K$, in Proposition \ref{AnalStrongDim3PGL2} we prove that $G$ is isomorphic to an open subgroup of a form of $\PGL_2$ over~$\BQ_p$. 
(An example for this is the group $H$ from Bartholdi-Bush \cite{BartholdiBush2007}; in retrospect it is thus no accident that $H$ was constructed within $\SL_2(\BQ_3)$.) 
Using this, the unramified Fontaine--Mazur conjecture then implies that any $p$-tower group $G_K$ with $G_K/D_4(G_K)$ in one of those 10 isomorphism classes is finite. 

%%%%%%%%%%%%%%%%%%%%%%%%
\medskip
%{\bf Heuristics for $G_K/D_4(G_K)$:} 
The second goal of this paper concerns the analogue of the Cohen--Lenstra heuristic for weak Schur $\sigma$-groups.
%for the quotients $G_K/D_4(G_K)$. 

\medskip
Recall that for any finite abelian $p$-group~$A$, Cohen and Lenstra \cite{CohenLenstra1984} predicted the frequency with which $A$ occurs as the maximal abelian quotient of $G_K$ among all imaginary quadratic fields~$K$. Boston, Bush, and Hajir \cite{BostonBushHajir2017} generalized this heuristic to the distribution of $G_K/P_j(G_K)$, where $P_j(G_K)$ denotes the $j$-th step of the lower $p$-central series of~$G_K$. In \cite{PinkRubio2025}, the second author and Rubio formalized this further by constructing a probability distribution on the set $\Sch$ of isomorphism classes of weak Schur $\sigma$-groups. 

This allows us to define the expected frequency for more general properties of~$G_K$, such as whether $G_K$ is finite or not. For instance, in the case $p>3$ it is proved in \cite[Thm.\,8.5]{Pink2025} that the set of isomorphism classes of infinite weak Schur $\sigma$-groups of Zassenhaus type (3,3) is a closed subset of $\Sch$ of measure~$0$. If the heuristic is correct, it thus follows that the ``if'' part of McLeman's statement holds on average, namely that almost all groups $G_K$ of Zassenhaus type $(3,3)$ are finite in a statistical sense, if $p>3$.
%the proportion of finite groups $G_K$ among those of Zassenhaus type $(3,3)$ tends to~$1$.
The results in this article now have the same consequences for $p=3$, provided that $G_K/D_4(G_K)$ lies in one of the 10 isomorphism classes above.

%%%%%%%%%%%%%%%%%%%%%%%%
\medskip
The formalism of \cite{PinkRubio2025} also yields an expected frequency for the isomorphism class of $G_K/E(G_K)$ for any canonically defined open normal subgroup $E(G_K)$. As in the original Cohen--Lenstra heuristic and in \cite{BostonBushHajir2017}, for a given finite $\sigma$-group $H$ this frequency is inversely proportional to the order of the group of $\sigma$-equivariant automorphisms of~$H$. 

\medskip
We were able to test this prediction experimentally in the case $E=D_4$ and $p=3$. This is based on an explicit presentation of the group $G_K/D_4(G_K)$ in terms of triple Massey products of cohomology classes in $H^1(G_K,\BZ/p\BZ)$ that was given by Vogel \cite[Prop.\,1.3.3]{Vogel-Massey}. After translating that into \'etale cohomology $H^1(\Spec\CO_K,\BZ/p\BZ)$ in Proposition \ref{prop:massey-presentation}, we can apply the explicit formulas for triple Massey products given by the first author and Carlson in~\cite{AhlqvistCarlson2025}. In the special case $d(G_K)=2$ the computation becomes particularly efficient, as explained in Section~\ref{Present}. Using the program \cite[{\tt 2-Massey}]{Ahlqvist-Pink-I}, we were thus able to compute the isomorphism class of $G_K/D_4(G_K)$ for all imaginary quadratic fields~$K$ with $d(G_K)=2$ and discriminant $-10^8 < d_K < 0$. 
These results match perfectly with the computations of the index-$p$-abel\-ian\-iza\-tion-data (IPAD) of Boston--Bush--Hajir \cite[Table~1]{BostonBushHajir2017}, 
%in the same range of discriminants $|d_K|<10^8$ as ours, 
as explained in Section~\ref{Exp}. 

In the case $p=3$ there are 19 possibilities for the isomorphism class of $G_K/D_4(G_K)$ with $d(G_K)=2$, including the 13 of Zassenhaus type $(3,3)$ that were mentioned above. For all 19 cases we compare the expected frequency with the observed frequency in Table \eqref{Experi} and find a reasonably good agreement. As the frequencies of the different cases vary over several orders of magnitude, we view this outcome as a moderate confirmation of the heuristic. 

The least likely of the 19 cases is that where $G_K/D_4(G_K)$ is isomorphic to $F_2/D_4(F_2)$, that is, where $G_K$ has Zassenhaus type $(a,b)$ with $a,b\ge5$. In this case $G_K$ is known to be infinite by Koch--Venkov \cite{KochVenkov1974}. Our computation shows that this case occurs for precisely 46 imaginary quadratic fields $K$ with $|d_K|<10^8$. The first examples of this and hence of infinite $G_K$ with $d(G_K)=2$ had already been found by the first author together with Carlson in \cite[Thm.\,5.8]{AhlqvistCarlson2025}.

%%%%%%%%%%%%%%%%%%%%%%%%
\medskip
In the future we plan to study some of the cases in more detail. For instance, GAP tells us already that a weak Schur $\sigma$-group $G$ with $p=3$ and $d(G)=2$ is isomorphic to $G/D_4(G)$ if and only if $G/D_4(G)$ is of type $[243,5]$ or $[243,7]$. The computation thus identifies all $83353+41398 = 124751$ imaginary quadratic fields $K$ with $|d_K| < 10^8$, whose $3$-tower group $G_K$ is finite with $d(G_K)=2$ and $|G_K|=3^5$.

%%%%%%%%%%%%%%%%%%%%%%%%
\medskip
%{\bf Content of the paper:}
The article is structured as follows. In Section~\ref{ProPGroups} we review mostly well-known material on $p$-groups and pro-$p$-groups. In Section~\ref{Powerful} we establish the criterion for an open normal subgroup of a finitely generated pro-$p$-group $G$ to be powerful. In Section \ref{SigmaAndSchur} we summarize basic properties of $\sigma$-pro-$p$-groups and of weak Schur $\sigma$-groups, taken mostly from \cite{PinkRubio2025}. In Section~\ref{AnalSchur} we prove a sufficient condition for a pro-$p$-group to be isomorphic to an open subgroup of a form of $\PGL_2$ over~$\BQ_p$ and deduce some consequences thereof. In Section~\ref{FinQuot} we discuss how to construct the finite quotients $G/E(G)$ of all weak Schur $\sigma$-groups of a certain type that are needed to apply the powerfulness criterion from Section~\ref{Powerful}. 
We expect this method and the corresponding code in the worksheet~\cite[{\tt 1-Checking-powerfulness}]{Ahlqvist-Pink-I} to be useful for future studies of weak Schur $\sigma$-groups and their finite quotients.
%the finite groups $G/P_j(G)$
In Section~\ref{Present} we describe an explicit presentation of $G_K/D_4(G_K)$ via triple Massey products in the case that the cup products of any two elements of $H^1(G_K,\BZ/p\BZ)$ vanish.
%and sketch the algorithm for its computation.

All sections up to this point concern an arbitrary odd prime~$p$, but in the remaining sections we assume that $p=3$. In Section~\ref{33ByD4} we analyze the different possibilities of $G/D_4(G)$ for all weak Schur $\sigma$-groups $G$ with $d(G)=2$, obtaining $19$ different isomorphism classes, of which $13$ have Zassenhaus type $(3,3)$. In Section~\ref{FrattPowerful} we explain the computer algebra computation which shows that $D_2(G)$ is powerful in $10$ of those $13$ cases. In Section~\ref{pTG} we outline the consequences of this for $p$-tower groups of imaginary quadratic fields. In the final Section~\ref{Exp} we present the experimental results concerning the isomorphism class of $G_K/D_4(G_K)$ for all imaginary quadratic fields $K$ with $|d_K| < 10^8$ and $d(G_K)=2$.

%%%%%%%%%%%%%%%%%%%%%%%%%%%%%%%%%%%%%%%%%%%%%%%%
%\newpage
\section{Pro-$p$-groups}
\label{ProPGroups}

In this section we fix an arbitrary prime $p$ and review some known facts about pro-$p$-groups.

%%%%%%%%%%%%%%%%%%%%%%%%
\medskip
Throughout this article, all homomorphisms of pro-$p$-groups are tacitly assumed to be continuous, and all subgroups are tacitly assumed to be closed. Thus by the subgroup generated by a subset of a pro-$p$-group $G$ we always mean the closure of the abstract subgroup generated by that set. In particular, when we say that a pro-$p$-group is generated by certain elements, we mean that it is topologically generated by them. For any subsets $A,B\subset G$ we let $[A,B]$ denote the subgroup generated by the commutators $[a,b]$ for all $a\in A$ and $b\in B$, which are defined by the convention $[a,b] := a^{-1}b^{-1}ab$. 
% Vogel uses the convention $[a,b] := a^{-1}b^{-1}ab$.}
We will need several canonically defined subgroups of~$G$. 

%%%%%%%%%%%%%%%%%%%%%%%%
\medskip
{\bf The Zassenhaus filtration:} 
%The \emph{lower central series of~$G$} is defined by $K_1(G) := G$ and $K_{i+1}(G) := [G,K_i(G)]$ for all $i\ge1$. (see Berkovich \cite[page xi]{Berkovich2008}). For any $i\ge0$ the $i$-th \emph{Agemo subgroup of~$G$} is the subgroup $\Agemo_i(G)$ generated by the subset $\{a^{p^i}\mid a\in G\}$. (see \cite[page xiii]{Berkovich2008}). All these are normal subgroups of~$G$.
For any $i\ge0$ the $i$-th \emph{Agemo subgroup of~$G$} is the normal subgroup $\Agemo_i(G)$ generated by the subset $\{a^{p^i}\mid a\in G\}$. 
The \emph{Zassenhaus filtration of~$G$} is a descending sequence of normal subgroups consisting of the \emph{dimension subgroups} $D_i(G)$ for all $i\ge1$ which, following a formula of Jennings \cite[Thm.\,12.9]{DdSMS2003}, can be characterized by the recursion relations
\UseTheoremCounterForNextEquation
\begin{equation}\label{Zass}
D_1(G) := G \quad\hbox{and}\quad 
D_i(G) := \Agemo_1\bigl(D_{\lceil i/p\rceil}(G)\bigr)\prod_{j+k=i}\bigl[D_j(G),D_k(G)\bigr]
 \quad\hbox{ for all}\ i\ge2.\kern-20pt
\end{equation}
%or by \cite[Thm.\,11.2]{DdSMS2003} as
% Formula of Lazard:
%$$D_i(G)\ =\ \prod_{jp^k\ge i} \Agemo_k(K_j(G)).$$
Each subquotient $\gr_i^D(G) := D_i(G)/D_{i+1}(G)$ is abelian and annihilated by~$p$ and hence an $\BF_p$-vector space. The subgroup $D_2(G)$ coincides with the Frattini subgroup $\Fratt(G)$, and the dimension $d(G)$ of $\gr^D_1(G)$ is the minimal number of generators of~$G$.
%We should avoid calling this the rank of~$G$, because the rank is sometimes defined as $\limsup_H d(H)$ as $H$ ranges over all open subgroups of~$G$.
The formation of $D_i(G)$ commutes with taking quotients, that is, for any normal subgroup $N\triangleleft G$ there is a natural isomorphism $D_i(G/N) \cong D_i(G)N/N$.

If $G$ is finitely generated, the $D_i(G)$ form a cofinal system of open normal subgroups of~$G$. In particular each $D_i(G)$ is then again finitely generated, 
%for instance by \cite[Prop.\,1.7]{DdSMS2003}
and each $\gr_i^D(G)$ has finite dimension over~$\BF_p$.

Also, the recursion relation \eqref{Zass} implies that $\bigoplus_{i\ge1}\gr_i^D(F_n)$ is a restricted Lie algebra over $\BF_p$ in the sense of \cite[\S12.1]{DdSMS2003} that is generated by $\gr_1^D(F_n)$. In particular this implies:

\begin{Prop}\label{AutFn}
Any automorphism of $G$ that acts trivially on $\gr_1^D(G)$ also acts trivially on $\gr_i^D(G)$ for every $i\ge1$.
\end{Prop}

%\begin{Proof}
%% I searched a little for a reference of this but eventually gave up.
%Let $I$ denote the augmentation ideal of the group ring $\BF_p[G]$. Then each $D_i(G)$ is defined as the subgroup of all $g\in G$ satisfying $g\equiv1$ modulo~$I^i$. In particular the automorphism in question acts trivially on $I/I^2$. The natural surjection $(I/I^2)^{\otimes i} \onto I^i/I^{i+1}$ thus implies that the automorphism also acts trivially on $I^i/I^{i+1}$. It therefore acts trivially on $\gr_i^D(F_n))$, as desired.
%\end{Proof}

%%%%%%%%%%%%%%%%%%%%%%%%
\medskip
{\bf The lower $p$-central series%
\footnote{Here we follow the same indexing as in O'Brien \cite[\S2]{OBrien1990}
%and Boston-Bush-Hajir \cite[\S2.1]{BostonBushHajir2017}
and the ANUPQ package of GAP,
%See Section 2.1-2 of the ANUPQ manual.
from which other sources such as Bartholdi-Bush \cite[p.\,161]{BartholdiBush2007} differ by a shift by~$1$.
% see also Berkovich \cite[\S26, bottom of p.264]{Berkovich2008}
%Section 39.17-13 of the GAP Reference Manual also has a shift by 1.
}:}
The recursion relations
\UseTheoremCounterForNextEquation
\begin{equation}\label{PjGDef}
P_0(G) := G \quad\hbox{and}\quad P_{j+1}(G) := \Agemo_1(P_j(G))[G,P_j(G)]
 \quad\hbox{ for all}\ j\ge0
\end{equation}
define another descending sequence of normal subgroups of~$G$, called the \emph{lower $p$-central series of~$G$.} The subgroup $P_1(G)$ coincides with the Frattini subgroup $\Fratt(G)$. For any normal subgroup $N\triangleleft G$ there is a natural isomorphism $P_j(G/N) \cong P_j(G)N/N$.

If $G$ is finitely generated, the $P_j(G)$ form a cofinal system of open normal subgroups of~$G$. If $G$ is finite, there is therefore a unique smallest integer $j$ with $P_j(G)=1$, which is called the \emph{$p$-class of~$G$}. Thus the trivial group has $p$-class~$0$, and a non-trivial elementary abelian $p$-group has $p$-class~$1$.
%Fortunately it seems that different sources agree on this, without a shift by 1.

%%%%%%%%%%%%%%%%%%%%%%%%
\medskip
{\bf Free pro-$p$-groups:}
For any integer $n\ge0$ we let $F_n$ denote the free pro-$p$-group on $n$ generators. Then a pro-$p$-group $G$ can be generated by $n$ elements if and only if there exists a surjective homomorphism $\phi\colon F_n\onto G$. This homomorphism is unique up to composition with an automorphism of~$F_n$, for instance by the same argument as in the proof of  \cite[Prop.\,6.3]{PinkRubio2025}.
% Koch \cite[Ch.\,4]{Koch2002} has a chapter on free pro-$p$-groups, but not what we want in a well citable form.
% Same for Neukirch-Schmidt-Wingberg \cite[\S3.5]{NSW2008}.
Thus we have $G\cong F_n/N$ for some normal subgroup $N\triangleleft F_n$, and $N$ is uniquely determined by $G$ up to the action of $\Aut(F_n)$. 

%\begin{Prop}\label{AutFn}
%The action of $\Aut(F_n)$ induces a surjective homomorphism 
%$$\Aut(F_n)\ \longonto\ \Aut(\gr_1^D(F_n))\ \cong\ \GL_n(\BF_p),$$
%and for any $i\ge1$ the action of $\Aut(F_n)$ on $\gr_i^D(F_n))$ factors through this homomorphism.
%\end{Prop}
%
%\begin{Proof}
%Part (a) follows from the fact that for any tuple of $n$ elements of $F_n$ which are linearly independent modulo $D_2(F_2)$, the group $F_n$ is again a free pro-$p$-group with these generators.
%REFERENCE?
%
%For (b) let $I$ denote the augmentation ideal of the group ring $\BF_p[G]$. Then each $D_i(G)$ is defined as the subgroup of all $g\in G$ with $g\equiv1$ modulo~$I^i$. The kernel of the homomorphism $\Aut(F_n)\onto \Aut(\gr_1^D(F_n))$ thus consists of all automorphisms which act trivially on $I/I^2$. By the natural surjection $(I/I^2)^{\otimes i} \onto I^i/I^{i+1}$, all these automorphisms also act trivially on $I^i/I^{i+1}$. They therefore also act trivially on $\gr_i^D(F_n))$, proving (b). 
%REFERENCE?
%\end{Proof}

%%%%%%%%%%%%%%%%%%%%%%%%
\medskip
{\bf Canonical subgroups:}
Consider a characteristic subgroup $E$ of~$F_n$, that is, a subgroup that is invariant under all automorphisms of~$F_n$. For any pro-$p$-group $G$ with $d(G)\le n$ choose a surjective homomorphism $\phi\colon F_n\onto\nobreak G$ with kernel~$N$. Then $E(G) := \phi(E)$ is a subgroup of~$G$ that is isomorphic to $EN/N$. Since $\phi$ is unique up to $\Aut(F_n)$, and $E$ is invariant under $\Aut(F_n)$, it follows that $E(G)$ is a characteristic subgroup of $G$ that depends only on $E$ and~$G$. The construction directly implies that the formation of $E(G)$ commutes with taking quotients, that is, for every normal subgroup $H\triangleleft G$ there is a natural isomorphism $E(G/H) \cong E(G)H/H$.

When $E$ is a proper subgroup of~$F_n$, the fact that $\Aut(F_n)$ surjects to $\Aut(F_n/\Fratt(F_n)) \cong \GL_n(\BF_p)$ implies that $E$ is contained in the Frattini subgroup $\Fratt(F_n)$. This then implies that $E(G)$ is contained in the Frattini subgroup $\Fratt(G)$, and so $G/E(G)$ has the same minimal number of generators as~$G$.

When $E$ is open in~$F_n$, the subgroup $E(G)$ is open in~$G$. This therefore provides a uniform way of constructing canonical open normal subgroups of $G$ with given properties. For example, if $E$ is any step in the Zassenhaus filtration of $F_n$, then $E(G)$ is the corresponding step in the Zassenhaus filtration of~$G$, and the same holds for the lower $p$-central series. But $E$ could also be constructed by any other combination of $p^i$-th powers and commutators, as in \eqref{E2Def} below.

%%%%%%%%%%%%%%%%%%%%%%%%
\medskip
{\bf The relative Frattini subgroup:} 
For any normal subgroup $N\triangleleft F_n$ we set
\UseTheoremCounterForNextEquation
\begin{equation}\label{NstarDef}
N^*\ :=\ \Agemo_1(N)[F_n,N].
\end{equation}
This is again a normal subgroup of $F_n$ and contained in~$N$. 
%In fact, it is the intersection of $N$ with all maximal proper subgroups of $N$ that are normal in~$F_n$.
Note that the recursion relation in \eqref{PjGDef} amounts to $P_{j+1}(F_n) = P_j(F_n)^*$ for all $j\ge0$.

\begin{Prop}\label{MstarProperties}
\begin{enumerate}
\item[(a)] For any normal subgroups $M,N\triangleleft F_n$ we have $(MN)^*=M^*N^*$.
\item[(b)] For any normal subgroups $M\subset N$ of~$F_n$, we have $M=N$ if and only if $MN^*=N$. 
\end{enumerate}
\end{Prop}

\begin{Proof}
In (a) the inclusion $\supset$ is obvious. To prove the other inclusion observe that $M^*N^*$ is a normal subgroup of~$F_n$. Thus every element of the quotient $F_n/M^*N^*$ commutes with the images of $M$ and $N$ and hence with $MN$, proving that $[F_n,MN] \subset M^*N^*$. In particular, any elements $m\in M$ and $n\in N$ commute with each other modulo $M^*N^*$; hence $(mn)^p$ is congruent to $m^pn^p$ modulo $M^*N^*$ and therefore lies in $M^*N^*$. Together this shows that $(MN)^*\subset M^*N^*$, as desired.

For (b) see for instance Neukirch-Schmidt-Wingberg \cite[Cor.\,3.9.3]{NSW2008}.
% also \cite[Lem.\,12.1.1]{WilsonProfiniteGroups}
\end{Proof}

\begin{Prop}\label{EstarGAndEG}
Consider a pro-$p$-group $G$ with $d(G)\le n$ and a characteristic subgroup $E$ of~$F_n$. Then we have $E(G)=1$ if and only if $E^*(G)=E(G)$.
\end{Prop}

\begin{Proof}
Write $G\cong F_n/N$ for a normal subgroup $N\triangleleft F_n$. Then by Proposition \ref{MstarProperties} (a) we have $N(EN)^* = NE^*N^* = E^*N$. By Proposition \ref{MstarProperties} (b) this implies that $N=EN$ if and only if $E^*N=EN$. Taking the respective quotients by~$N$, the proposition follows.
\end{Proof}

%%%%%%%%%%%%%%%%%%%%%%%%
\medskip
{\bf The minimal number of relations:} 
Consider a pro-$p$-group $G$ with $d(G)=n$, and choose a normal subgroup $N\triangleleft F_n$ with $F_n/N\cong G$. Then the construction \eqref{NstarDef} implies that $N/N^*$ is abelian and annihilated by~$p$ and hence an $\BF_p$-vector space. Proposition \ref{MstarProperties} (b) implies that its dimension is the minimal number of generators of $N$ as a normal subgroup of~$F_n$. In other words $r(G) := \dim_{\BF_p}(N/N^*)$ is \emph{the minimal number of relations defining $G$ as a quotient of~$F_n$}. 
%In [Koch 1997 Algebraic Number Theory, Ch.3 \S1.13] this is called the \emph{relation rank of~$G$}.

As a variant of this, consider a characteristic subgroup $E$ of~$F_n$. Then $E(G)=1$ is equivalent to $E\subset N$; and in that case we can form the $\BF_p$-vector space 
\UseTheoremCounterForNextEquation
\begin{equation}\label{RelEGDef}
\Rel_E(G)\ :=\ N/EN^*.
\end{equation}
Its dimension $r_E(G)$ depends only on $E$ and $G$ and is \emph{the minimal number of relations defining $G$ as a quotient of~$F_n/E$}. 
%Though the construction of $\Rel_E(G)$ depends on the chosen presentation of~$G$, this should not cause any problems below.
%In Wilson \cite[beginning of \S12.1]{WilsonProfiniteGroups} this would be $d_{F_n/E}(N/E)$.

\begin{Prop}\label{RelRel}
For any proper characteristic subgroups $E\subset D$ of~$F_n$ and any pro-$p$-group $G$ with $d(G)=n$ we have
$$r_E(G/E(G))\ \ge\ r_D(G/D(G)).$$
\end{Prop}

\begin{Proof}
Choose a normal subgroup $N\triangleleft F_n$ such that $G\cong F_n/N$. Then the assumption $d(G)=n$ implies that $N$ is contained in the Frattini subgroup $\Fratt(F_n)$. Since the proper characteristic subgroups $E$ and $D$ are also contained in $\Fratt(F_n)$, the quotients $G/E(G) \cong F_n/EN$ and $G/D(G) \cong F_n/DN$ again have minimal number of generators~$n$.

Now Proposition \ref{MstarProperties} (a) implies that $(EN)^*=E^*N^*$ and hence $E(EN)^*=EN^*$. From \eqref{RelEGDef} we therefore obtain $\Rel_E(G/E(G)) = EN/EN^*$. The same argument with $D$ in place of~$E$ yields $\Rel_D(G/D(G)) = DN/DN^*$. We thus have a natural surjective homomorphism
\UseTheoremCounterForNextEquation
\begin{equation}\label{RelRelIsom}
\Rel_E(G/E(G))\ =\ EN/EN^*\ \longonto\ DN/DN^*\ =\ \Rel_D(G/D(G)).
\end{equation}
This directly implies the desired inequality.
\end{Proof}

%%%%%%%%%%%%%%%%%%%%%%%%
\medskip
{\bf The $p$-cover:} 
Now let $G$ be a finite $p$-group with $d(G)=n$. Then $G\cong F_n/N$ for an open normal subgroup $N\triangleleft F_n$. The normal subgroup $N^*$ defined in \eqref{NstarDef} is then again open, and the finite $p$-group $G^* := F_n/N^*$ is called the \emph{$p$-cover of~$G$}. By construction the kernel of the projection $G^*\onto G$ is isomorphic to $N/N^*$ and therefore a finite dimensional $\BF_p$-vector space contained in the center of~$G^*$. The extension $G^*\onto G$ possesses a universal property with respect to all central extensions of $G$ by $\BF_p$-vector spaces. In particular, its isomorphism class is therefore independent of the choice of the presentation of~$G$. 

Following O'Brien \cite[\S2]{OBrien1990}, the $\BF_p$-vector space $N/N^*$ is called the \emph{$p$-multiplicator of~$G$}, and its dimension $\mu(G)$ is called the \emph{$p$-multiplicator rank of~$G$}. This coincides with the minimal number $r(G)$ of relations of $G$ as a quotient of~$F_n$, but we will keep both notations because they are relevant in different contexts. 
%Actually, they appear in the same kind of formula (d$'$) and (d$''$).
Also, if $G$ has $p$-class~$j$, then $P_j(G^*)$ is contained in $\Ker(G^*\,{\onto}\,G) \cong N/N^*$ and hence an $\BF_p$-subspace. This is called the \emph{nucleus of~$G$} and its dimension is called the \emph{nuclear rank of~$G$}. 
% O'Brien only defines the multiplicator and the nucleus at the top of page 680, and the multiplicator rank only sort of at the beginning of \S3.1. He never defines the nuclear rank, but uses the term at the beginning of \S5.

%\begin{Prop}\label{PGdim}
%Let $G$ be a finite $p$-group of $p$-class~$\le j$. Then $P_j(G^*)$ is a finite dimensional $\BF_p$-vector space of dimension
%$$\dim_{\BF_p} P_j(G^*)\ =\ 
%\scriptstyle\biggl\{\displaystyle\begin{array}{cl}
%\nu(G) &\hbox{if $G$ has $p$-class~$j$,}\\[3pt]
%0 &\hbox{if $G$ has $p$-class $<j$.}
%\end{array}$$
%\end{Prop}
%
%\begin{Proof}
%If $G$ has $p$-class~$j$, this is just the definition of $\nu(G)$. Otherwise we have $P_{j-1}(G)=\nobreak 1$. Writing $G\cong F_n/N$ as above, this means that $P_{j-1}(F_n)\subset M$. But this implies that $P_j(F_n) = P_{j-1}(F_n)^* \subset M^*$ and hence $P_j(G^*) = P_j(F_n)M^*/M^*=1$. Thus its dimension is~$0$.
%\end{Proof}

\begin{Prop}\label{MuRelE}
For any proper characteristic subgroup $E$ of~$F_n$ and any pro-$p$-group $G$ with $d(G)=n$ and $E(G)=1$ the group $E(G^*)$  is a finite dimen\-sional $\BF_p$-vector space with
$$\mu(G)\ =\ r_E(G) + \dim_{\BF_p}E(G^*).$$
\end{Prop}

\begin{Proof}
(Compare Boston-Nover \cite[Prop.\,2]{BostonNover2006}.) 
With $G\cong F_n/N$ as above, the assumption $E(G)=1$ is equivalent to $E\subset N$. We thus have the inclusions $N^* \subset EN^* \subset N$. Since $G$ is finite, the quotient $N/N^*$ is a finite dimensional $\BF_p$-vector space; hence the same follows for all subquotients of these inclusions. Counting their dimensions, this implies that
$$\dim_{\BF_p}(N/N^*)\ =\ \dim_{\BF_p}(N/EN^*) + \dim_{\BF_p}(EN^*/N^*).$$
Here by definition the left hand side is $\mu(G)$, the first term on the right hand side is $r_E(G)$, and the equality $G^*=F_n/N^*$ implies that $E(G^*) = EN^*/N^*$.
\end{Proof}

%%%%%%%%%%%%%%%%%%%%%%%%
\medskip
{\bf Descendants:} Consider now a finite $p$-group~$G$ of $p$-class $j>0$. Then a finite $p$-group $H$ such that $H/P_j(H) \cong G$ is called a \emph{descendant of~$G$}. If in addition $H$ has $p$-class $j+1$, then $H$ is called an \emph{immediate descendant of~$G$}. All immediate descendants are quotients of $G^*$ by a subgroup of $\Ker(G^*{\,\onto\,}G)$. O'Brien's $p$-group generation algorithm \cite{OBrien1990} provides an efficient way of constructing all descendants of $G$ of a given type, up to isomorphism.

%%%%%%%%%%%%%%%%%%%%%%%%%%%%%%%%%%%%%%%%%%%%%%%%
%\newpage
\section{Criterion for an open normal subgroup to be power\-ful}
\label{Powerful}

{}From now on we assume that the prime $p$ is odd. Consider a finitely generated pro-$p$-group~$G$, and recall that its Frattini subgroup is $\Fratt(G) = \Agemo_1(G)[G,G]$. Following Berkovich \cite[\S26]{Berkovich2008} the group $G$ is called \emph{powerful} if and only if $\Agemo_1(G) = \Fratt(G)$.
%$K_2(G) \subset \Agemo_1(G)$.
This notion is relevant for us because, by a theorem of Lazard \cite{Lazard1965}, \cite[Thm.\,8.1]{DdSMS2003}, a topological group has the structure of a $p$-adic analytic group if and only if it possesses an open subgroup which is a powerful finitely generated pro-$p$-group.
% by the following theorem of Lazard \cite{Lazard1965}, see also Dixon-du Sautoy-Mann-Segal \cite[Thm.\,8.1]{DdSMS2003}:
%
%\begin{Thm}\label{PowerfulAnalytic}
%A topological group has the structure of a $p$-adic analytic group if and only if it possesses an open subgroup which is a powerful finitely generated pro-$p$-group.
%\end{Thm}

%%%%%%%%%%%%%%%%%%%%%%%%
\medskip
In this section we give a criterion for an open normal subgroup of $G$ to be powerful.

\begin{Prop}\label{HisPowerful}
For any open normal subgroup $H\triangleleft G$, the quotient
$$\bar H\ :=\ H/\Agemo_1(H)[G,\Fratt(H)]$$
is finite, and $H$ is powerful if and only if $\Fratt(\bar H)=1$.
%That is also equivalent to $\bar H$ being powerful.
\end{Prop}

\begin{Proof}
Since $H$ is open in the finitely generated pro-$p$-group~$G$, it is itself finitely generated. Thus its Frattini subgroup $\Fratt(H)$ is open in~$H$. Repeating the argument shows that $\Fratt(\Fratt(H))$ is open in $\Fratt(H)$ and hence in~$H$. Since $\Fratt(\Fratt(H)) =  \Agemo_1(\Fratt(H))[\Fratt(H),\Fratt(H)])$ is contained in $\Agemo_1(H)[G,\Fratt(H)]$, the latter subgroup is thus also open in~$H$. The quotient $\bar H$ is therefore finite.

Now by definition $H$ is powerful if and only if the quotient $\Fratt(H)/\Agemo_1(H)$ is trivial. Since a $p$-group is trivial if and only if its Frattini quotient is trivial, that is equivalent to $\Fratt(H)/\Agemo_1(H)\Fratt(\Fratt(H))$ being trivial. That is a finite dimensional $\BF_p$-vector space with a continuous action of the pro-$p$-group~$G$, because $H$ is normal in~$G$. Thus that quotient is trivial if and only if its group of $G$-coinvariants is trivial. This group of coinvariants is naturally isomorphic to 
$$\Fratt(H)\!\bigm/\!\Agemo_1(H)\Fratt(\Fratt(H))[G,\Fratt(H)] 
\ \cong\ \Fratt(H)\!\bigm/\!\Agemo_1(H)[G,\Fratt(H)].$$
But that is just the Frattini subgroup $\Fratt(\bar H)$ of~$\bar H$, and we are done.
\end{Proof}

%%%%%%%%%%%%%%%%%%%%%%%%
\medskip
Now we fix an integer $n\ge0$ and an open characteristic subgroup $E$ of~$F_n$. To this we associate the subgroups
\UseTheoremCounterForNextEquation
\begin{eqnarray}\label{E1Def}
E_1 &\!:=\!& \Fratt(E), \\[3pt]
\UseTheoremCounterForNextEquation\label{E2Def}
E_2 &\!:=\!& \Agemo_1(E)[F_n,\Fratt(E)].
\end{eqnarray}
By construction these are again characteristic subgroups of~$F_n$. Applying Proposition \ref{HisPowerful} to $(G,H)=(F_n,E)$ shows that $E_2$ is open in~$F_n$. Since $E_1$ contains~$E_2$, it is also open. By Section~\ref{ProPGroups}, to any pro-$p$-group $G$ with $d(G)\le n$, the subgroups $E_2\subset E_1\subset E$ therefore associate canonical open normal subgroups $E_2(G)\subset E_1(G)\subset E(G)$ of~$G$.

\begin{Prop}\label{EGisPowerful}
For any pro-$p$-group $G$ with $d(G)\le n$, the subgroup $E(G)$ is powerful if and only if $E_2(G)=E_1(G)$.
\end{Prop}

\begin{Proof}
Choose a normal subgroup $N\triangleleft F_n$ and an isomorphism $G\cong F_n/N$. Consider the subgroup $H:=E(G)\cong EN/N$. Then by \eqref{E2Def} its quotient from Proposition \ref{HisPowerful} is 
$$\bar H\ :=\ H/\Agemo_1(H)[G,\Fratt(H)]
\ \cong\ EN/\Agemo_1(E)[F_n,\Fratt(E)]N
\ =\ EN/E_2N.$$
Therefore
$$\Fratt(\bar H)\ \cong\ \Fratt(E)N/E_2N\ =\ E_1N/E_2N\ \cong\ E_1(G)/E_2(G).$$
The assertion is now a direct consequence of Proposition \ref{HisPowerful}.
\end{Proof}

%%%%%%%%%%%%%%%%%%%%%%%%
\medskip
Since the condition in Proposition \ref{EGisPowerful} depends only on the finite quotient $G/E_2(G)$, it is possible to test it with a computer algebra system for a large number of groups~$G$.

%%%%%%%%%%%%%%%%%%%%%%%%%%%%%%%%%%%%%%%%%%%%%%%%
%\newpage
\section{$\sigma$-Pro-$p$-groups and weak Schur $\sigma$-groups}
\label{SigmaAndSchur}

As before we fix an odd prime~$p$. In \cite{PinkRubio2025} the second author and Luca Rubio developed the concept of weak Schur $\sigma$-groups and their classifying space. To avoid repetition, we refer to \cite[\S2--4]{Pink2025} for the basic notation and results that we need here.

%%%%%%%%%%%%%%%%%%%%%%%%
\medskip
{\bf $\sigma$-Pro-$p$-groups:} As a brief reminder, a \emph{$\sigma$-pro-$p$-group} is a pro-$p$-group together with an action of a group $\{1,\sigma\}$ of order~$2$, and a finite $\sigma$-pro-$p$-group is called a \emph{$\sigma$-$p$-group}. A $\sigma$-equivariant homomorphism between $\sigma$-pro-$p$-groups is called a \emph{$\sigma$-homomorphism}. The notions of \emph{$\sigma$-isomorphisms} and \emph{$\sigma$-automorphisms} are defined accordingly. The group of $\sigma$-automorphisms of a $\sigma$-pro-$p$-group $G$ is denoted $\Aut_\sigma(G)$.

To any $\sigma$-pro-$p$-group $G$ we associate the natural closed subsets
\UseTheoremCounterForNextEquation
\begin{eqnarray}\label{G+Def}
G^+ &\!:=\!& \bigl\{ a\in G\bigm| {}^\sigma a = a\bigr\}, \\[3pt]
\UseTheoremCounterForNextEquation\label{G-Def}
G^- &\!:=\!& \bigl\{ a\in G\bigm| {}^\sigma a = a^{-1}\bigr\}.
\end{eqnarray}
We call the elements of $G^+$ \emph{even} and the elements of $G^-$ \emph{odd}. We call $G$ \emph{totally even} if $G=G^+$, and \emph{totally odd} if $G=G^-$. If $G$ is abelian, these are subgroups that form a $\BZ/2\BZ$-grading of~$G$. In the general case one still has a weaker kind of grading by \cite[Prop.\,3.4]{Pink2025}. Also the formation of $G^\pm$ is exact in short exact sequences by \cite[Prop.\,3.3]{Pink2025}.

%%%%%%%%%%%%%%%%%%%%%%%%
\medskip
For every integer $n\ge0$ we turn the free pro-$p$-group on $n$ generators $F_n$ into a $\sigma$-pro-$p$-group by making the generators odd. Thus any $\sigma$-pro-$p$-group that is generated by $n$ odd elements is $\sigma$-isomorphic to  $F_n/N$ for a $\sigma$-invariant normal subgroup $N\triangleleft F_n$. For these groups we have the following useful fact:

\begin{Prop}\label{IsomVSSigmaIsom}
For any $\sigma$-pro-$p$-group $G$ that is generated by finitely many odd elements, its $\sigma$-iso\-mor\-phism class is determined uniquely by its isomorphism class as a pro-$p$-group.
\end{Prop}

\begin{Proof}
By the argument at the bottom of page 640 of \cite{BostonBushHajir2017}, all automorphisms of order $2$ of~$G$ that act as $-1$ on the abelianization $G_\ab$ are conjugate in $\Aut(G)$.
\end{Proof}

\medskip
Also, any characteristic subgroup $E$ of~$F_n$ is automatically $\sigma$-invariant. Thus for any $\sigma$-pro-$p$-group $G$ that is generated by $n$ odd elements, the associated characteristic subgroup $E(G)$ of $G$ is automatically $\sigma$-invariant. In the framework of $\sigma$-pro-$p$-groups, the construction of $E(G)$ actually only requires that $E$ is a normal subgroup of $F_n$ that is invariant under $\Aut_\sigma(F_n)$, but in practice that will not make a difference for us.

%%%%%%%%%%%%%%%%%%%%%%%%
%\medskip
%{\bf Canonical subgroups:}
%For any $\sigma$-invariant normal subgroup $E$ of~$F_n$ that is invariant under $\Aut_\sigma(F_n)$, we can repeat the construction of canonical subgroups from Section~\ref{ProPGroups}. Namely, for any $\sigma$-pro-$p$-group $G$ that is generated by $n$ odd elements, choose a surjective $\sigma$-homomorphism $\phi\colon F_n\onto G$ with kernel~$N$. Then $E(G) := \phi(E) \cong EN/N$ is a $\sigma$-invariant normal subgroup of~$G$. Since $\phi$ is unique up to $\Aut_\sigma(F_n)$ by \cite[Prop.\,6.3]{PinkRubio2025}, and $E$ is invariant under $\Aut_\sigma(F_n)$, it follows that $E(G)$ depends only on $E$ and~$G$. Again this construction commutes with taking quotients, and in the case of an open subgroup $E$ of~$F_n$ it yields an open subgroup $E(G)$ of~$G$.

%%%%%%%%%%%%%%%%%%%%%%%%
\medskip
{\bf Schur $\sigma$-groups:} A \emph{weak Schur $\sigma$-group} is a $\sigma$-pro-$p$-group $G$ with the properties
\UseTheoremCounterForNextEquation
\begin{equation}\label{SchConsDef}
\scriptstyle\biggl\{\displaystyle
\begin{array}{c}
\hbox{$\dim_{\BF_p}\!H^2(G,\BF_p) \le \dim_{\BF_p}\!H^1(G,\BF_p) < \infty$} \\[3pt]
\hbox{and $\sigma$ acts by $-1$ on both spaces}
\end{array}\scriptstyle\biggr\}.
\end{equation}
The weak Schur $\sigma$-groups $G$ with $\dim_{\BF_p}\!H^1(G,\BF_p) \le n$ are precisely those that are $\sigma$-isomorphic to $F_n/N$ for a $\sigma$-invariant normal subgroup $N\triangleleft F_n$ that is generated by all conjugates of $n$ odd elements of~$F_n$. In \cite[\S7--9]{PinkRubio2025} the set $\Sch$ (pronounce \emph{Schur}) of $\sigma$-isomorphism classes of weak Schur $\sigma$-groups was endowed with a natural topology and a probability measure~$\mu_\infty$. 
More precisely, for every integer $n\ge0$ the subset 
\UseTheoremCounterForNextEquation
\begin{equation}\label{Schn}
\Sch_n\ :=\ \bigl\{[G]\in\Sch \bigm| d(G)=n\bigr\}
\end{equation}
is open and closed in~$\Sch$. Moreover, the subsets \cite[(10.1)]{PinkRubio2025}
\UseTheoremCounterForNextEquation
\begin{equation}\label{UDG}
U_{D,G}\ :=\ \bigl\{[H]\in\Sch_n \bigm| H/D(H)\cong G/D(G) \bigr\}
\end{equation}
for all open proper characteristic subgroups $D$ of~$F_n$ and all $[G]\in\Sch_n$ are open and closed and form a basis for the topology of $\Sch_n$. Furthermore, for any integer $k\ge0$ as well as for $k=\infty$ consider the positive real number 
\UseTheoremCounterForNextEquation
\begin{equation}\label{CkDef}
C_k\ :=\ \prod_{i=1}^k (1-p^{-i}).
\end{equation}
Then by \cite[Prop.\,10.2]{PinkRubio2025} we have 
\UseTheoremCounterForNextEquation
\begin{equation}\label{MuUDG}
\mu_\infty(U_{D,G})\ =\ \frac{C_\infty}{C_{n-m}}\cdot\frac{1}{|\Aut_\sigma(G/D(G))|}\,,
\end{equation}
where $m=r_D(G/D(G))$ is the minimal number of relations defining $G/D(G)$ as a quotient of~$F_n/D$. In particular, by \cite[Prop.\,9.8]{PinkRubio2025} we have 
\UseTheoremCounterForNextEquation
\begin{equation}\label{MuSchn}
\mu_\infty(\Sch_n)\ =\ \frac{C_\infty}{C_n^2} \cdot p^{-n^2}.
\end{equation}

Following \cite[Def.\,11.1]{PinkRubio2025}, a weak Schur $\sigma$-group $G$ such that for every open subgroup $H<G$ the abelianization $H_\ab$ is finite is called a \emph{strong Schur $\sigma$-group}. If one only requires that $G_\ab$ is finite one gets a \emph{Schur $\sigma$-group} by \cite[Prop.\,11.2]{PinkRubio2025}. Clearly every finite weak Schur $\sigma$-group is a strong Schur $\sigma$-group. In \cite[Thm.\,11.4]{PinkRubio2025} we proved that the subset $\Sch^\strong$ of $\sigma$-isomorphism classes of strong Schur $\sigma$-groups has measure~$1$.

%%%%%%%%%%%%%%%%%%%%%%%%%%%%%%%%%%%%%%%%%%%%%%%%
%\newpage
\section{Analytic strong Schur $\sigma$-groups}
\label{AnalSchur}

As before we assume that $p$ is odd. It seems an interesting general problem to classify all infinite weak Schur $\sigma$-groups that have the structure of a $p$-adic analytic group. Here we deal only with a very special case.

\begin{Prop}\label{AnalStrongDim3PGL2}
Let $G$ be an infinite pro-$p$-group with $d(G)=2$ and the properties:
\begin{enumerate}
\item[(a)] Its Frattini subgroup $\Fratt(G)$ is powerful with $d(\Fratt(G))\le3$.
\item[(b)] For every open subgroup $H$ of $G$ the abelianization $H_\ab$ is finite.
\end{enumerate}
Then $G$ is isomorphic to an open subgroup of a form of $\PGL_2$ over~$\BQ_p$.
\end{Prop}

\begin{Proof}
Abbreviate $N:=\Fratt(G)$. Then condition (a) and \cite[Thm.\,3.8]{DdSMS2003} imply that every open subgroup $H$ of $N$ satisfies $d(H) \le d:= d(N)\le 3$. Combining \cite[Thm.\,4.2, Def.\,4.7, Thm.\,8.36]{DdSMS2003} it follows that $N$ and hence $G$ is a $p$-adic analytic group of dimension $\le3$. Since $G$ is infinite, its dimension is also $\ge1$. 

Let $\Fg$ denote the $\BQ_p$-Lie algebra associated to~$G$. Then there exist an open subgroup $L\subset G$ and a $\BZ_p$-lattice $\Fl\subset\Fg$ such that the exponential map induces a homeomorphism $\Fl\isoto L$. This map sends $\Fl\cap[\Fg,\Fg]$ to a normal subgroup $K\triangleleft L$ and induces a group isomorphism $\Fl/(\Fl\cap[\Fg,\Fg]) \isoto L/K$. In particular $L/K$ is abelian; hence it must be finite by the second assumption on~$G$. This proves that $[\Fg,\Fg]=\Fg$.

But we already know that $1\le \dim_{\BQ_p}(\Fg)\le 3$. Thus $\Fg$ must have a simple semisimple quotient of dimension $\le3$. By the classification of semisimple Lie algebras this must be a form of ${\Fs\Fl}_2$ over~$\BQ_p$ and have dimension $=3$. It follows that $d=3$ and that $\Fg$ itself is a form of ${\Fs\Fl}_2$ over~$\BQ_p$.

Its automorphism group $\CG := \Aut(\Fg)$ is therefore a form of $\PGL_2$ over~$\BQ_p$. Moreover the adjoint representation $\Ad_G\colon G\to\CG(\BQ_p)$ induces an isomorphism from $L$ to an open subgroup of $\CG(\BQ_p)$. Thus the kernel $\Ker(\Ad_G)$ is finite and its image is an open subgroup of $\CG(\BQ_p)$. It remains to show that $\Ker(\Ad_G)=1$. 

For this we recall that, since $N$ is powerful, by \cite[Lem.\,26.8~(a)]{Berkovich2008} we have $P_{i+1}(N) = \Fratt(P_i(N))$ for every $i\ge0$%
\footnote{Berkovich \cite[\S26, bottom of p.264]{Berkovich2008} also uses a different convention for the $p$-central series.}.
Thus $\gr_i^P(N) := P_i(N)/P_{i+1}(N)$ is the Frattini quotient of $P_i(N)$ and hence an $\BF_p$-vector space of dimension $d(P_i(N))$. Since $N$ is a $p$-adic analytic group of dimension $d=3$, this dimension is equal to $3$ for every $i\gg0$. On the other hand, again since $N$ is powerful, by \cite[Lem.\,26.8 (b)]{Berkovich2008} the map $x\mapsto x^p$ induces a surjective homomorphism $\gr_i^P(N) \onto \gr_{i+1}^P(N)$ for every $i\ge0$. Since $d(N)\le3$ by assumption, it follows that this map is an isomorphism for every $i\ge0$. 

Now suppose that $\Ker(\Ad_G)$ is non-trivial, and choose an element $x$ of order~$p$. If this does not lie in $N$, the fact that $G$ has minimal number of generators $d(G)=2$ implies that $G$ is generated by $x$ and one other element $y\in G\setminus N$. But then $\Ad_G(G)$ is generated by the image of $y$ alone and can therefore not be open in $\CG(\BQ_p)$. Thus $x$ must lie in~$N$. 

Since the subgroups $P_i(N)$ for all $i\ge0$ form a neighborhood base of the identity in~$N$, there then exists an $i\ge0$ with $x\in P_i(N)\setminus P_{i+1}(N)$. Thus  $x$ represents a non-zero element of $\gr_i^P(N)$. Since $x\mapsto x^p$ induces an isomorphism $\gr_i^P(N) \onto \gr_{i+1}^P(N)$, it follows that $x^p\not=1$, contradicting the assumption on~$x$. Thus $\Ker(\Ad_G)$ is trivial, and we are done.
\end{Proof}

%%%%%%%%%%%%%%%%%%%%%%%%

\begin{Prop}\label{AnalStrongDim3Countable}
There exist at most countably many isomorphism classes of groups $G$ as in Proposition~\ref{AnalStrongDim3PGL2}.
\end{Prop}

\begin{Proof}
Every form $\CG$ of $\PGL_2$ over $\BQ_p$ arises from a quaternion algebra $Q$ over~$\BQ_p$ such that $\CG(\BQ_p) \cong Q^\times/\BQ_p^\times$. As there are precisely two isomorphism classes of quaternion algebras over a local field, 
%for instance by \cite[Thm.\,13.1.6]{Voight2021}. 
there are only two isomorphism classes for~$\CG$. Since the topology of $\CG(\BQ_p)$ is second countable, there exist at most countably many open compact subgroups.
\end{Proof}

%%%%%%%%%%%%%%%%%%%%%%%%
\medskip
From this we can deduce that infinite weak Schur $\sigma$-groups with the property \ref{AnalStrongDim3PGL2} (a) form a sparse subset of~$\Sch$:

\begin{Prop}\label{AnalDim3WeakSchurNullSet}
The set of $\sigma$-isomorphism classes of infinite weak Schur $\sigma$-groups $G$ with $d(G)=2$, whose Frattini subgroup $\Fratt(G)$ is powerful with $d(\Fratt(G))\le3$, is a closed subset of $\Sch$ of measure~$0$.
\end{Prop}

\begin{Proof}
First, the set of $[G]\in\Sch$ with $d(G)=2$ is open and closed by the construction of the topology on $\Sch$. Next set $E:=\Fratt(F_2)$ and define open characteristic subgroups $E_2\subset E_1:=\Fratt(E)$ of~$F_n$ as in \ref{E1Def} and \ref{E2Def}. For any weak Schur $\sigma$-group $G$ with $d(G)=2$ we then have open characteristic subgroups $E_2(G)\subset E_1(G) \subset E(G)=\Fratt(G)$ of~$G$. 
By Proposition \ref{EGisPowerful} the condition that $\Fratt(G)$ is powerful depends only on the isomorphism class of $G/E_2(G)$. Moreover, since $E_1(G)=\Fratt(\Fratt(G))$, the number $d(\Fratt(G)) = \dim_{\BF_p} \Fratt(G)/\Fratt(\Fratt(G))$ also depends only on the isomorphism class of $G/E_2(G)$.
Since for any finite group~$H$, the set of $\sigma$-isomorphism classes of weak Schur $\sigma$-groups $G$ with $d(G)=2$ and $G/E_2(G)\cong H$ is open and closed in $\Sch$ by \cite[(10.1)]{PinkRubio2025}, it follows that the set of $\sigma$-isomorphism classes of weak Schur $\sigma$-groups $G$ with $d(G)=2$, whose Frattini subgroup $\Fratt(G)$ is powerful with $d(\Fratt(G))\le3$,
is open and closed in~$\Sch$.
On the other hand, the subset of $\sigma$-isomorphism classes of all finite Schur $\sigma$-groups is open in $\Sch$ by \cite[Prop.\,8.7]{PinkRubio2025}. Together this implies that the set in question is a closed subset of~$\Sch$.

To show that it has measure~$0$ recall that the set of $\sigma$-isomorphism classes of weak Schur $\sigma$-groups that are not strong has measure~$0$ by  \cite[Thm.\,11.4]{PinkRubio2025}. Thus it suffices to show that the set of $\sigma$-isomorphism classes of strong Schur $\sigma$-groups in our set has measure~$0$. By Proposition \ref{AnalStrongDim3Countable} there are at most countably many isomorphism classes of such groups. By Proposition \ref{IsomVSSigmaIsom} the same thus also holds for their $\sigma$-isomorphism classes. Moreover, for every infinite Schur $\sigma$-group $G$ the subset $\{[G]\}\subset\Sch$ has measure $0$ by \cite[Prop.\,4.7]{Pink2025}. As the measure is countably additive, it follows that the our set has measure~$0$, as desired.
\end{Proof}

\section{Finite quotients of weak Schur $\sigma$-groups}
\label{FinQuot}

In this section we fix an integer $n\ge0$ and concentrate on weak Schur $\sigma$-groups $G$ with minimal number of generators $d(G)=n$. For any open proper characteristic subgroup $E$ of~$F_n$, we discuss how to generate all possibilities for $G/E(G)$, using recursion on~$E$. For this the following abbreviation will be helpful:

\begin{Def}\label{SchurEQuotDef}
A finite $p$-group $H$ is called a \emph{Schur $E$-quotient} if and only if there exists a weak Schur $\sigma$-group $G$ with $d(G)=n$ such that $G/E(G) \cong H$. 
\end{Def}

Note that, since $E$ is contained in $\Fratt(F_n)$, any Schur $E$-quotient $H$ has minimal number of generators $d(H)=n$. 

%%%%%%%%%%%%%%%%%%%%%%%%
\medskip
For the recursion we postulate that we already have a list of possibilities for some initial choice of~$E$. For the induction step we therefore take two open proper characteristic subgroups $E\subset D$ of~$F_n$ and suppose that we have a list of Schur $D$-quotients $H$ up to isomorphism. Any Schur $E$-quotient $K\cong G/E(G)$ then satisfies $K/D(K)\cong H$ for the Schur $D$-quotient $H:=G/D(G)$. For any $H$ in our list we must therefore find all Schur $E$-quotients $K$ up to isomorphism with $K/D(K)\cong H$. In principle O'Brien's $p$-group generation algorithm is well-suited for a problem of this kind. However, we know of no implementation that deals with $\sigma$-$p$-groups.

Fortunately, since any Schur $E$-quotient is generated by odd elements, by Proposition \ref{IsomVSSigmaIsom} its $\sigma$-iso\-mor\-phism class is determined uniquely by its isomorphism class as a $p$-group without~$\sigma$. Instead of carrying along the action of~$\sigma$, we can therefore work with $p$-groups alone, testing independently whether they possess \emph{some} structure of a $\sigma$-$p$-group that is generated by odd elements. 

But we must also guarantee that the desired Schur $E$-quotient~$K$ can be described by odd relations as a quotient of $F_n/E$. We have not investigated how to decide this efficiently for a given finite $\sigma$-$p$-group. Instead, we restrict ourselves to the case that the minimal number of relations of $H$ as a quotient of $F_n/D$ is already the maximal permitted value of~$n$ for Schur $D$-quotients. Then, if a covering group $K$ has the same minimal number of relations as a quotient of $F_n/E$, using \eqref{RelRelIsom} we find that those relations can also be taken as odd. With this strategy we arrive at the following formal criterion:

\begin{Prop}\label{HKGEDProp}
Take open proper characteristic subgroups $E\subset D$ of $F_n$. Consider a Schur $D$-quotient~$H$ with $r_{D}(H)= n$ and a $p$-group~$K$ with $d(K)=n$ and $K/D(K)\cong\nobreak H$. Then $K$ is a Schur $E$-quotient if and only if the following conditions hold: 
\begin{enumerate} 
\item[(a)] We have $E(K)=1$.
\item[(b)] The group $K$ possesses the structure of a $\sigma$-$p$-group that is generated by odd elements.
\item[(c)] We have $r_E(K)=n$.
\end{enumerate}
\end{Prop}

\begin{Proof}
First we show:

\begin{Lem}\label{HKGEDLem1}
The ``only if'' part of Proposition \ref{HKGEDProp} holds.
\end{Lem}

\begin{Proof}
If $K$ is a Schur $E$-quotient, take a weak Schur $\sigma$-group $G$ with $d(G)=n$ and $G/E(G) \cong K$. Then (a) and (b) follow automatically. 
Next we claim that
$$n\ \ge\ r(G)\ \ge\ r_E(K)\ \ge\ r_D(H)\ =\ n.$$
Indeed, the first inequality follows from the definition of weak Schur $\sigma$-groups. The second inequality follows by applying Proposition \ref{RelRel} with $(G,1,E)$ in place of $(G,E,D)$. Also the isomorphy $K/D(K) \cong H$ implies that $G/D(G) \cong\nobreak H$; hence the third inequality holds by Proposition \ref{RelRel}. Finally, the last equality holds by assumption. Combining the inequalities implies that $r_E(K)=n$, proving~(c).
\end{Proof}

\medskip
For the rest of the proof we assume conditions (a--c). Fixing a structure of $\sigma$-$p$-group on $K$ as in condition (b), there is then a unique structure of $\sigma$-$p$-group on $H$ for which the given isomorphism $K/D(K)\cong H$ is $\sigma$-equivariant. Since $d(K)=n$, we can choose a $\sigma$-invariant open normal subgroup $L\triangleleft F_n$ and a $\sigma$-isomorphism $K\cong F_n/L$. By (a) this subgroup then satisfies $E\subset L$. Moreover, the $\sigma$-isomorphism $K/D(K)\cong H$ yields a $\sigma$-isomorphism $H\cong F_n/M$ for the $\sigma$-invariant open normal subgroup $M:=DL$. 

\begin{Lem}\label{HKGEDLem2}
The $\sigma$-$p$-group $\Rel_D(H)$ is a totally odd $\BF_p$-vector space of dimension~$n$.
\end{Lem}

\begin{Proof}
By assumption there exists a weak Schur $\sigma$-group $G$ with $d(G)=n$ such that $G/D(G) \cong H$. Both sides of this isomorphism are $\sigma$-$p$-groups that are generated by odd elements, so by Proposition \ref{IsomVSSigmaIsom} we may assume the isomorphism to be $\sigma$-equivariant. Since $d(G)=\nobreak n$, any $\sigma$-homomorphism $F_n\to G$ lifting the $\sigma$-epimorphism $F_n\onto F_n/M\cong H$ is then automatically surjective. Its kernel is therefore a $\sigma$-invariant normal subgroup $N\triangleleft\nobreak F_n$ with a $\sigma$-isomorphism $G\cong F_n/N$, such that $DN=\nobreak M$. The map \eqref{RelRelIsom} from the proof of Proposition \ref{RelRel} with $(G,1,D)$ in place of $(G,E,D)$ is then a $\sigma$-epimorphism of $\BF_p$-vector spaces
$$\Rel_1(G)\ =\ N/N^*\ \longonto\ M/DM^*\ =\ \Rel_D(H).$$
But by the definition of weak Schur $\sigma$-groups, the ideal $N$ is generated by odd relations; hence $N/N^*$ is totally odd. The same therefore follows for the right hand side. As its dimension is $r_D(H)= n$ by assumption, we are done.
\end{Proof}

\begin{Lem}\label{HKGEDLem3}
The ``if'' part of Proposition \ref{HKGEDProp} holds.
\end{Lem}

\begin{Proof}
Since $E(K)=1$, we have $K/E(K)\cong K$. The map \eqref{RelRelIsom} from the proof of Proposition \ref{RelRel} with $(K,E,D)$ in place of $(G,E,D)$ is thus an surjective homomorphism of $\BF_p$-vector spaces
$$\Rel_E(K)\ =\ L/EL^*\ \longonto\ M/DM^*\ =\ \Rel_D(H).$$
Here the right hand side has dimension $r_D(H)=n$ by the assumption on~$H$, and the left hand side has dimension $r_E(K)=n$ by condition~(c). Thus the homomorphism is an isomorphism. But by construction it is also $\sigma$-equivariant; hence it is a $\sigma$-isomorphism. By Lemma \ref{HKGEDLem2} it thus follows that $L/EL^*$ is a totally odd $\BF_p$-vector space of dimension~$n$.

Now recall that the formation of the odd part is exact in short exact sequences of $\sigma$-pro-$p$-groups by \cite[Prop.\,3.3]{Pink2025}. We can therefore lift $n$ generators of $L/EL^*$ to $n$ odd elements of~$F_n$. Letting $N$ be the subgroup generated by their conjugates, the quotient $G:=F_n/N$ is then a weak Schur $\sigma$-group with $d(G)\le n$. Moreover, by construction we have $N\subset L$ and $ENL^* = L$, and by Proposition \ref{MstarProperties} (b) this implies that $EN = L$. This means that $G/E(G) \cong F_n/EN = F_n/L \cong K$. Finally, since $d(K)=n$, this isomorphy implies that $d(G)=n$; hence $K$ is a Schur $E$-quotient, as desired.
\end{Proof}

\medskip
Combining Lemmas \ref{HKGEDLem1} and \ref{HKGEDLem3} now finishes the proof of Proposition \ref{HKGEDProp}.
\end{Proof}

%%%%%%%%%%%%%%%%%%%%%%%%
\medskip
For any given subgroups $E\subset D$ as in Proposition \ref{HKGEDProp}, we can choose open characteristic subgroups $E=E_m\subset E_{m-1}\subset\ldots\subset E_0=D$ of $F_n$ satisfying $E_{j-1}^*\subset E_j$ for all $1\le j\le m$. To generate all possibilities for $G/E(G)$ with given $G/D(G)$, it then suffices to apply the proposition inductively with $E_j\subset E_{j-1}$ in place of $E\subset D$ for all~$j$.

We are thus interested in Proposition \ref{HKGEDProp} under the additional assumption $D^*\subset E$. In that case $\Ker(K\!\onto\!H)$ is a finite dimensional $\BF_p$-vector space and we can replace condition \ref{HKGEDProp}~(c) by the computationally more efficient condition \ref{HKGDstarEDProp} (c$'$) below. Its advantage lies in the fact that $\mu(H)$ and $\dim_{\BF_p} E(H^*)$ depend only on the group~$H$, so the condition gives an a priori value for the step size $\dim_{\BF_p}\Ker(K\!\onto\!H)$. 

\begin{Prop}\label{HKGDstarEDProp}
Take open proper characteristic subgroups $E\subset D$ of $F_n$ with $D^*\subset\nobreak E$. Consider a Schur $D$-quotient~$H$ with $r_{D}(H)= n$ and a $p$-group~$K$ with $d(K)=n$ and $K/D(K)\cong\nobreak H$. Then $K$ is a Schur $E$-quotient if and only if the following conditions hold: 
\begin{enumerate} 
\item[(a)] We have $E(K)=1$.
\item[(b)] The group $K$ possesses the structure of a $\sigma$-$p$-group that is generated by odd elements.
\item[(c$\kern2pt'$)] We have \ $\mu(H) = \dim_{\BF_p}\Ker(K\!\onto\!H) + n + \dim_{\BF_p} E(H^*).$
\end{enumerate}
\end{Prop}

\begin{Proof}
If $K$ satisfies (a), as in the proof of Proposition \ref{HKGEDProp} we can choose an open normal subgroup $L\triangleleft F_n$ such that $K\cong F_n/L$. This subgroup  then satisfies $E\subset L$, and the isomorphy $K/D(K)\cong H$ implies that $H\cong F_n/M$ for the open normal subgroup $M:=DL$. 
By Proposition \ref{MstarProperties}~(a) this last equality implies that $M^* = D^*L^*$. Using $D^*\subset E$ this shows that $EM^* = ED^*L^* = EL^* \subset L$. Together we therefore have the inclusions 
$$M^*\ \subset\ EM^*\ =\ EL^*\ \subset\ L\ \subset\ M.$$
Since $H\cong F_n/M$ is finite, the quotient $M/M^*$ is a finite dimensional $\BF_p$-vector space. The same thus follows for all subquotients of the above inclusions. Counting their dimensions, this implies that
%\UseTheoremCounterForNextEquation
%\begin{equation}\label{HKGDstarEDPropEq1}
$$\dim_{\BF_p}(M/M^*)\ =\  \dim_{\BF_p}(M/L) + \dim_{\BF_p}\bigl(L/EL^*\bigr) + \dim_{\BF_p}(EM^*/M^*).$$
%\end{equation}
We claim that we can rewrite this equation in the form
\UseTheoremCounterForNextEquation
\begin{equation}\label{HKGDstarEDPropEq2}
\mu(H) \ =\ \dim_{\BF_p}\Ker(K\!\onto\!H) + r_E(K) + \dim_{\BF_p} E(H^*).
\end{equation}
Indeed, the left hand sides agree by the definition of the $p$-multiplicator rank. Next, the short exact sequence $1\to M/L\to F_n/L\to F_n/M\to 1$ implies that $M/L \cong\hbox{$\Ker(K\!\onto\!H)$}$; hence first terms on the right hand sides are equal. The following two terms coincide by the definition \eqref{RelEGDef} of $\Rel_E(K)$. Finally, the defining equality $F_n/M^* = H^*$ of the $p$-cover implies that $EM^*/M^* = E(H^*)$, so the last terms match as well.

In particular this proves that all terms in the equation \eqref{HKGDstarEDPropEq2} are finite. 
The equation thus implies that (c$'$) is equivalent to \ref{HKGEDProp}~(c), making the proposition a special case of Proposition~\ref{HKGEDProp}. 
\end{Proof}

\medskip
One can use Proposition~\ref{HKGDstarEDProp} to construct all Schur $P_j(F_n)$-quotients by induction on~$j$, and in that case the condition (c$'$) simplifies further, as follows.
%This is useful, because the $p$-multiplicator rank and the nuclear rank can be computed efficiently.

\begin{Prop}\label{HKGDstarPPProp}
Take any integer $j\ge1$. Consider a Schur $P_j(F_n)$-quotient~$H$ with $r_{P_j(F_n)}(H) = n$ and a $p$-group~$K$ with $d(K)=n$ and $K/P_j(K)\cong\nobreak H$. Then $K$ is a Schur $P_{j+1}(F_n)$-quotient if and only if the following conditions hold: 
\begin{enumerate} 
\item[(a)] We have $P_{j+1}(K)=1$.
\item[(b)] The group $K$ possesses the structure of a $\sigma$-$p$-group that is generated by odd elements.
\item[(c$\kern2pt''$)] We have \ $\mu(H) = \dim_{\BF_p}\Ker(K\!\onto\!H) + n.$
\end{enumerate}
\end{Prop}

\begin{Proof}
Take $D:=P_j(F_n)$ and $E:=P_{j+1}(F_n)$. Then $H$ has $p$-class $\le j$, and so its $p$-cover $H^*$ has $p$-class $\le j+1$;
%Writing $H\cong F_n/M$, this means that $P_j(F_n)\subset M$. But that implies that $P_{j+1}(F_n) = P_j(F_n)^* \subset M^*$ and hence $P_{j+1}(H^*) = P_{j+1}(F_n)M^*/M^*=1$. 
hence we have $E(H^*)=P_{j+1}(H^*)=1$. The condition \ref{HKGDstarEDProp}~(c$'$) thus reduces to condition~(c$''$). The result is therefore a special case of Proposition \ref{HKGDstarEDProp}.
\end{Proof}

\begin{Rem}\label{HKGDstarPPNuProp}\rm
When $H$ has $p$-class~$j$, by Proposition \ref{MuRelE} we have $\mu(H) = r_{P_j(F_n)}(H) + \dim_{\BF_p}P_j(H^*) = n + \nu(H)$; hence (c$''$) is further equivalent to $\dim_{\BF_p}\Ker(K\!\onto\!H) = \nu(H)$.
\end{Rem}

%CAUTION: When $H$ has $p$-class $j'<j$, then $\nu(H) = \dim_{\BF_p}P_{j'}(H^*)$ and the formula changes.

%NOTE FOR US: I checked with ANUPQ that the $p$-multiplicator rank and the nuclear rank are computed efficiently for about 500 groups of $p$-class $j\le6$. The computation via computing the $p$-cover explicitly took more than $10^4$ times as long as via the inbuilt functions {\tt MultiplicatorRank} and {\tt NuclearRank}.

%%%%%%%%%%%%%%%%%%%%%%%%%%%%%%%%%%%%%%%%%%%%%%%%
%\newpage
\section{An explicit presentation via Massey products}
\label{Present}

As before we fix an odd prime~$p$. Consider a number field $K$ and let $\widetilde K$ denote its maximal unramified pro-$p$-extension. Following McLeman \cite{McLeman2008} the group $G_K := \Gal(\widetilde K/K)$ is called the \emph{$p$-tower group} associated to~$K$. In this section we use the formulas for triple Massey products from \cite{AhlqvistCarlson2025} to derive an explicit presentation of $G_K/D_4(G_K)$ in the case when all cup products of any two elements of $H^1(G_K,\BZ/p\BZ)$ vanish.

\medskip
First recall that the maximal abelian quotient of $G_K$ is naturally isomorphic to the $p$-primary part of the ideal class group $\Cl(K)$; hence the minimal number of generators $d(G_K)$ is equal to $\dim_{\BF_p}{\Cl}(K)[p]$. Next, by Koch \cite[Thm.\,7.23]{Koch2002}, the group $G_K/P_2(G_K)$ has a presentation in terms of the cup product 
\UseTheoremCounterForNextEquation
\begin{equation}\label{Bockstein}
\smallsmile\,:\ H^1(G_K, \BZ/p\BZ) \times H^1(G_K, \BZ/p\BZ) \longto H^2(G_K,\BZ/p\BZ)
\end{equation}
and the Bockstein homomorphism 
%\UseTheoremCounterForNextEquation
%\begin{equation}\label{Bockstein}
$B\colon H^1(G_K, \BZ/p\BZ)\to H^2(G_K, \BZ/p\BZ).$
%\end{equation}
%which is the connecting homomorphism in the long exact cohomology sequence associated to the short exact sequence $0\to\BZ/p\BZ\to\BZ/p^2\BZ\to\BZ/p\BZ\to0$.

In \cite[Prop.\,4.2]{AhlqvistCarlsonOpen} the first author and Carlson showed how to compute the cup product in terms of \'etale cohomology. When all cup products vanish, this implies in particular that $G_K/D_3(G)$ is isomorphic to $F_n/D_3(F_n)$ for $n:= d(G_K)$. 

As explained in \cite[Section~2]{AhlqvistCarlson2025} or \cite{DwyerMassey}, for any classes $x,y,z\in H^1(G_K, \BZ/p\BZ)$ whose cup products $x\smallsmile y$ and $y\smallsmile z$ both vanish, one can define the triple \emph{Massey product} $\langle x,y,z\rangle$, which in general is a subset of $H^2(G_K, \BZ/p\BZ)$. If the cup products of \emph{all} elements of $H^1(G_K, \BZ/p\BZ)$ vanish, one can therefore define the Massey product of any triple of elements in $H^1(G_K, \BZ/p\BZ)$, and in that case it is actually single-valued and defines a trilinear map
\UseTheoremCounterForNextEquation
\begin{equation}\label{Massey}
\langle\ ,\ ,\ \rangle:\ H^1(G_K, \BZ/p\BZ) \times H^1(G_K, \BZ/p\BZ) \times H^1(G_K, \BZ/p\BZ) \longto H^2(G_K,\BZ/p\BZ).
\end{equation}
Under this assumption there is an analogous presentation of $G_K/D_4(G_K)$ in terms of 
%the Bockstein map and the 
triple Massey products due to Vogel \cite[Prop.\,1.3.3]{Vogel-Massey}. 

%%%%%%%%%%%%%%%%%%%%%%%%
\medskip
We will need to reformulate this presentation in terms of \'etale cohomology. For this write $X:=\Spec \sO_K$. Then we may identify $G_K$ with the maximal pro-$p$-quotient of the \'etale fundamental group of~$X$, or equivalently with the Galois group of the maximal unramified pro-$p$-cover $\widetilde X := \Spec\CO_{\widetilde K}$ of~$X$. By Milne \cite[Remark 2.21(b)]{MilneEtale} this cover determines a spectral sequence in \'etale cohomology
\[E_2^{i,j} = H^i(G_K, H^j(\widetilde X, \BZ/p\BZ)) \Longrightarrow H^{i+j}(X, \BZ/p\BZ),\]    
and hence a natural homomorphism
\UseTheoremCounterForNextEquation
\begin{equation}\label{HiGKHiX}
\kappa_i\colon H^i(G_K, \BZ/p\BZ) \to H^i(X,\BZ/p\BZ)
\end{equation}
for every $i\ge0$. Since $H^1(\widetilde X,\BZ/p\BZ)=0$ by construction, this homomorphism is an isomorphism when $i=1$ and an injection when $i=2$. 

The Bockstein homomorphism, the cup product, and the Massey product can all be formed with the same properties within the \'etale cohomology $H^*(X,\BZ/p\BZ)$, and each of these constructions is compatible with the homomorphisms \eqref{HiGKHiX}. In particular, the cup products of any two elements of $H^1(G_K,\BZ/p\BZ)$ vanishes if and only if the same holds for any two elements of $H^1(X,\BZ/p\BZ)$. For more details see \cite[Section~3]{AhlqvistCarlson2025}.

We will also need the \emph{Artin-Verdier duality} \cite{ArtinVerdierSeminar1964}, \cite[2.4]{MazurEtale}, which gives a perfect pairing 
\UseTheoremCounterForNextEquation
\begin{equation}\label{ArtinVerdier}
\langle\ ,\ \rangle:\ H^i(X, \BZ/p\BZ)\times H^{3-i}(X,D(\BZ/p\BZ)) \longto H^3(X, D(\BZ/p\BZ))\xrightarrow[\sim]{\operatorname{inv}}\BZ/p\BZ\,,
\end{equation}
of \'etale cohomology groups, where $D(\BZ/p\BZ):=R\HOM(\BZ/p\BZ, \mathbb{G}_m)$ is an object in the derived category $D(X_{\textsf{\'et}})$. 
% We may also identify $H^{j}(X,D(\BZ/p\BZ))$ with $H^j(X_{\textsf{Fppf}}, \pmb{\mu}_p)$ where we use the big fppf site, as in \cite[Remark 4.1]{AhlqvistCarlson2025}.

%%%%%%%%%%%%%%%%%%%%%%%%
\medskip
Now write $F_n$ as the free pro-$p$-group on $n$ generators $a_1,\ldots,a_n$. For any indices $i,j,k$ we then have $[[a_i,a_j],a_k]\in D_3(F_n)$, and since $p\ge3$, we also have $a_i^p\in D_3(F_n)$. Since $\gr_3^D(F_n)$ is an $\BF_p$-vector space, it follows that for integers $r$ the residue classes of $[[a_i,a_j],a_k]^r$ and $a_i^{pr}$ modulo $D_4(F_n)$ depend only on $r$ modulo~$p$. Thus the following result makes sense:

\begin{Prop}\label{prop:massey-presentation}
Assume that the cup product of any two elements in $H^1(G_K, \BZ/p\BZ)$ vanishes. Choose $\BF_p$-bases $x_1,\ldots,x_n$ of $H^1(X, \BZ/p\BZ)$ and $b_1,\ldots,b_m$ of $H^{1}(X,D(\BZ/p\BZ))$. Then $G_K/D_4(G_K)$ is isomorphic to $F_n/D_4(F_n)N$, where $N$ is the normal subgroup generated by the elements $f_\ell$ for $1\leq \ell\leq m$ that are defined by
$$f_\ell\ =
\prod_{1\leq i\leq n} \hskip-2pt a_i^{-pe_{iii}(\ell)} \hskip-5pt
\prod_{\scriptstyle1\leq i<j\leq n \atop\scriptstyle1\leq k\leq j} \hskip-5pt [[a_i,a_j],a_k]^{e_{ijk}(\ell)}\,,$$
with  
$$e_{ijk}(\ell)\ =\ \scriptstyle\biggl\{\displaystyle\hskip-3pt
\begin{array}{rl}
-\langle\langle x_k,x_i,x_j\rangle,b_\ell\rangle & \text{if } k<j, \\[3pt] 
\phantom{-}\langle\langle x_j,x_j,x_i\rangle,b_\ell\rangle & \text{if } k=j\,,
\end{array}$$
%$$f_\ell\ = 
%\prod_{1\leq i\leq n} \hskip-2pt a_i^{p\langle B(x_i), b_\ell\rangle}\hskip-5pt
%\prod_{\scriptstyle 1\leq i<j\leq n \atop \scriptstyle 1\leq k<j} \hskip-5pt 
%   [[a_i,a_j],a_k]^{-\langle\langle x_j,x_i,x_k\rangle,b_\ell\rangle}\hskip-5pt
%\prod_{1\leq i<j\leq n} \hskip-5pt 
%   [[a_i,a_j],a_j]^{\langle\langle x_i,x_j,x_j\rangle,b_\ell\rangle},\hskip-20pt$$
where the products are taken in any order and where by abuse of notation we replace the terms in the exponents by any representatives in~$\BZ$.
%Eric: By Vogel 1.2.15, the Bockstein is $B(x)=-\langle x,x,x\rangle$ when cup products vanish, so it is cleaner to only mention Massey products instead of Bockstein + Massey.
\end{Prop}

Note that for $p>3$ we have $a_i^p\in D_4(F_n)$ for all~$i$; hence in that case the first product in the formula for $f_\ell$ can be left out.

%%%%%%%%%%%%%%%%%%%%%%%%
\bigskip
\begin{Proof}
Let $\chi_1, \dots \chi_n$ be the basis of $H^1(G_K, \BZ/p\BZ)$ mapping to $x_1,\dots, x_n$ under the isomorphism $\kappa_1$ from~\eqref{HiGKHiX}. Then the compatibility of the Massey product with the maps \eqref{HiGKHiX} means that $\kappa_2(\langle\chi_i,\chi_j,\chi_k\rangle) = \langle x_i,x_j,x_k\rangle$ for all $i,j,k$.

Next, by dualizing, the Artin-Verdier duality \eqref{ArtinVerdier} and the injection $\kappa_2$ yield a surjection
$$\xymatrix@C-0pt{
H^1(X,D(\BZ/p\BZ)) \ar@{}[r]|-{\textstyle\cong} &
H^2(X,\BZ/p\BZ)^\vee \ar@{->>}[r]^-{\kappa_2^\vee} &
H^2(G_K,\BZ/p\BZ)^\vee. \\}$$
The given basis $b_1,\ldots,b_m$ of $H^{1}(X,D(\BZ/p\BZ))$ therefore yields generators $\phi_1,\ldots,\phi_m$ of $H^2(G_K,\BZ/p\BZ)^\vee$ which by construction satisfy $\langle\kappa_2(\xi),b_\ell\rangle = \phi_\ell(\xi)$ for all $\xi\in H^2(G_K,\BZ/p\BZ)$.
% Third, the compatibility of the respective Bockstein homomorphisms with the maps $\kappa_i$ implies that 
% $$\langle B(x_i),b_\ell\rangle\ =\  
% \langle B(\kappa_1(\chi_i))),b_\ell\rangle\ =\  
% \langle \kappa_2(B(\chi_i))),b_\ell\rangle\ =\  
% \phi_\ell(B(\chi_i)).$$

Third, choose elements $g_1,\ldots,g_n\in G_K$ whose residue classes modulo $\Fratt(G_K)$ form the dual basis of the basis $\chi_1, \dots, \chi_n$ of $H^1(G_K, \BZ/p\BZ) \cong (G_K/\Fratt(G_K))^\vee$. Then these elements generate $G_K$ and therefore define a presentation $1\to R\to F_n\to G_K\to 1$ with $a_i\mapsto g_i$ for all~$i$. From the Hochschild--Serre spectral sequence we then obtain the transgression isomorphism 
$$\tra\colon (R/R^*)^\vee\ \cong\ H^1(R, \BZ/p\BZ)^{G_K}\isoto H^2(G_K, \BZ/p\BZ)\,,$$
since $H^2(F_n, \BZ/p\BZ)=0$ by \cite[Thm.\,4.12]{Koch2002}. Its dual map $H^2(G_K,\BZ/p\BZ)^\vee \to R/R^*$ thus sends $\phi_1,\ldots,\phi_m$ to a system of generators  of $R/R^*$. Writing these as the residue classes of $f_1, \dots, f_m \in R$, it follows that $f_1,\ldots,f_m$ generate $R$ as a normal subgroup of~$F_n$. By Vogel \cite[Prop.\,1.3.3, Ex.\,1.2.11]{Vogel-Massey} the statement of the proposition therefore holds with
$$e_{ijk}(\ell)\ =\ \scriptstyle\biggl\{\displaystyle\hskip-3pt
\begin{array}{rcrl}
\varphi_\ell(-\langle \chi_k,\chi_i,\chi_j\rangle) 
&\kern-5pt=\kern-5pt& -\langle\langle x_k,x_i,x_j\rangle,b_\ell\rangle 
%\varphi_\ell(-\langle \chi_j,\chi_i,\chi_k\rangle) 
%&\kern-5pt=\kern-5pt& -\langle\langle x_j,x_i,x_k\rangle,b_\ell\rangle 
& \text{if } k<j, \\[3pt] 
\varphi_\ell(\langle \chi_j,\chi_j,\chi_i\rangle) 
&\kern-7pt=\kern-7pt& \langle\langle x_j,x_j,x_i\rangle,b_\ell\rangle 
%\varphi_\ell(\langle \chi_i,\chi_j,\chi_j\rangle) 
%&\kern-7pt=\kern-7pt& \langle\langle x_i,x_j,x_j\rangle,b_\ell\rangle 
& \text{if } k=j,
\end{array}$$
as desired.
\end{Proof}

%%%%%%%%%%%%%%%%%%%%%%%%
\medskip
In the special case $n:=d(G_K)=2$ the formula from Proposition \ref{prop:massey-presentation} reduces to 
\UseTheoremCounterForNextEquation
\begin{equation}\label{PresentationDim2}
f_\ell\ =\ 
a_1^{-p \langle \langle x_1,x_1,x_1\rangle, b_\ell\rangle}
 a_2^{-p \langle \langle x_2,x_2,x_2\rangle, b_\ell\rangle}
 [[a_1, a_2], a_1]^{-\langle\langle x_1, x_1, x_2\rangle, b_\ell\rangle}
 [[a_1, a_2], a_2]^{\langle\langle x_2, x_2, x_1\rangle, b_\ell\rangle}.\hskip-15pt
\end{equation}
This involves only Massey products of the form $\langle x,x,y\rangle$, which can be computed more efficiently than Massey products of three distinct elements, using the formulas in \cite{AhlqvistCarlson2025} and~\cite{AhlqvistCarlson2023}.

\medskip
Namely, nontrivial classes $x\in H^1(X, \BZ/p\BZ)$ are represented by pairs $(L_x,\rho_x)$ consisting of an unramified cyclic extension $L_x/K$ of degree $p$ and a generator $\rho_x$ of $\Gal(L/K)$. On the other hand, classes in $H^1(X, D(\BZ/p\BZ))$ can be represented by pairs $(a, J)$ consisting of a fractional ideal $J$ of $\CO_K$ and an element $a\in K^\times$ with $(a)=J^{-p}$.
%Another such pair $(a',J')$ represent the same class if and only if there exists $b\in K^\times$ with $(a',J')=(ab^{-p},bJ)$. 
By \cite[Thm.\,4.12]{AhlqvistCarlson2025} we can then compute $\langle\langle x,x,y\rangle, (a,J)\rangle$ using PARI as follows:
\begin{enumerate}
\item Find solutions $I$ to the equation $i_x(J)+(1-\rho_x)^2I=0$ in $\Cl(L_x)$ using {\tt matsolvemod}. 
\item Choose a solution $I$ such that the principal ideal $(i_x(J)+(1-\rho_x)^2I)^{-1}$ has a generator $t\in L_x^\times$ (using {\tt bnfisprincipal0}) whose norm is $N_{L_x/K}(t)=au^p$ for some unit $u$ in~$\CO_K^\times$ (using {\tt bnfisunit0}).  
\item Letting $\rho_y^{-1}\colon\Gal(L_y/K)\isoto\BZ/p\BZ$ denote the isomorphism with $\rho_y\mapsto[1]$ and $\operatorname{Art}_{L_y}\colon \Cl(K)\to \Gal(L_y/K)$ the Artin symbol, we then have\\[-25pt]
\end{enumerate}  
\UseTheoremCounterForNextEquation
\begin{equation}\label{MasseyProductFormula}
\langle\langle x,x,y\rangle, (a,J)\rangle=\rho_y^{-1}\bigl(\operatorname{Art}_{L_y}(N_{L_x/K}(I)+J)\bigr).
\end{equation}

%%%%%%%%%%%%%%%%%%%%%%%%
\medskip
All this applies to imaginary quadratic fields, because:

\begin{Prop}\label{prop:cup-vanish} % Also true for real quadratic fields
When $K$ is an imaginary quadratic field, 
\begin{enumerate} 
\item[(a)] the cup product of any two elements in $H^1(G_K, \BZ/p\BZ)$ vanishes, and
\item[(b)] the map $\kappa_i\colon H^i(G_K, \BZ/p\BZ) \to H^i(X,\BZ/p\BZ)$ is an isomorphism for $i=1,2$, unless $p=3$ and $K=\mathbb{Q}(\sqrt{-3})$.
\end{enumerate}
\end{Prop}

\begin{Proof}
In this case $G_K$ is a Schur $\sigma$-group by Koch-Venkov \cite[\S1]{KochVenkov1974}; hence $\sigma$ acts by $-1$ on $H^1(G_K,\BZ/p\BZ)$ and $H^2(G_K,\BZ/p\BZ)$ by \eqref{SchConsDef}. Thus $\sigma$ acts by $+1$ on $\bigwedge^2H^1(G_K,\BZ/p\BZ)$. As the cup product commutes with~$\sigma$, and the prime $p$ is odd, the cup product must therefore vanish completely. This proves (a). For (b), see \cite[\S5, p.\,100]{AhlqvistCarlson2025}. 
% Next we already know that $\kappa_i$ is an isomorphism for $i=1$ and injective for $i=2$. Moreover $H^2(X,\BZ/p\BZ)$ is dual to $H^1(X,D(\BZ/p\BZ))$, which lies in a short exact sequence $1\to \CO_K^\times/(\CO_K^\times)^p \to H^1(X,D(\BZ/p\BZ)) \to \Cl(K)[p] \to 1$ (see \cite[p.~14]{AhlqvistCarlson2025}).  
% Since $K$ is an imaginary quadratic field and $p$ is odd, the group $\CO_K^\times/(\CO_K^\times)^p$ is trivial except in the case $p=3$ and $K=\mathbb{Q}(\sqrt{-3})$. With that exception we therefore have 
% $$\dim_{\BF_p} H^2(X,\BZ/p\BZ) = \dim_{\BF_p} \Cl(K)[p] = \dim_{\BF_p} \Cl(K)/p\Cl(K) = \dim_{\BF_p} H^1(G_K, \BZ/p\BZ).$$
% But since $G_K$ is a Schur $\sigma$-group, this number is also equal to $\dim_{\BF_p} H^2(G_K, \BZ/p\BZ)$. 
% %\cite{Shafarevich1963}
% The injective homomorphism $\kappa_2$ is therefore an isomorphism, finishing the proof of~(b).
\end{Proof}

%%%%%%%%%%%%%%%%%%%%%%%%
\medskip
The program \cite[{\tt 2-Massey}]{Ahlqvist-Pink-I} does the necessary computations for any imaginary quadratic field $K$ with $d(G_K)=2$, using the PARI library \cite{PARI}.
Beginning with the discriminant $d_K$ and an odd prime~$p$, the program finds a basis $x_1,x_2$ of $H^1(X, \BZ/p\BZ)$ via the PARI function {\tt bnrclassfield} and a basis $b_1,b_2$ of $H^1(X,D(\BZ/p\BZ))$, and computes the eight values $\langle\langle x_i,x_i,x_j\rangle, b_\ell\rangle$ in~\eqref{PresentationDim2} using the formula~\eqref{MasseyProductFormula}. 
From this it produces a list of defining relations for $G_K/D_4(G_K)$ in GAP syntax.

%%%%%%%%%%%%%%%%%%%%%%%%%%%%%%%%%%%%%%%%%%%%%%%%
%\newpage
\section{Weak Schur $\sigma$-groups with $d(G)=2$ for $p=3$}
\label{33ByD4}

From now on we assume that $p=3$. Any weak Schur $\sigma$-group $G$ with $d(G)=2$ is the quotient of $F_2$ by the normal subgroup generated by two odd relations contained in $D_2(F_2)$. As the subquotient $\gr^D_2(F_2)$ is totally even, these relations must then already lie in $D_3(F_2)$. To see what this means recall from McLeman \cite[Thm.\,5]{McLeman2009} that, since $p=3$, we have
\UseTheoremCounterForNextEquation
\begin{equation}\label{F2Grad}
\dim_{\BF_3}\gr^D_i(F_2)\ =\ 
\scriptstyle\left\{\kern-2pt\textstyle 
\begin{array}{ll}
2 &\hbox{if \ $i=1$,}\\[3pt]
1 &\hbox{if \ $i=2$,}\\[3pt]
4 &\hbox{if \ $i=3$.}
\end{array}\right.
\end{equation}
Thus the relations generate a subspace of dimension $\le2$ of the $\BF_3$-vector space $\gr^D_3(F_2)$ of dimension~$4$, and there are many possibilities for that. If its dimension is $=2$, the group $G$ is called \emph{of Zassenhaus type (3,3)} (see Koch-Venkov \cite{KochVenkov1974}). 
% A weak Schur $\sigma$-group~$G$ has Zassenhaus type (3,3) if and only if it can be described by two generators modulo two relations, which lie in the subgroup $D_3(F_2)$ for the Zassenhaus filtration of~$F_2$ and whose residue classes in the subquotient $\gr^D_3(F_2)$ are $\BF_3$-linearly independent. 

%%%%%%%%%%%%%%%%%%%%%%%%
\medskip
Some basic facts about such groups are:

\begin{Prop}\label{HisGbyD4}
For any finite $3$-group $H$ the following are equivalent:
\begin{enumerate}
\item[(a)] $H$ is isomorphic to $G/D_4(G)$ for a weak Schur $\sigma$-group $G$ with $d(G)=2$.
\item[(b)] $\dim_{\BF_3}\gr_1^D(H)=2$ and $\dim_{\BF_3}\gr_2^D(H)=1$ and $\dim_{\BF_3}\gr_3^D(H)\ge2$ and $D_4(H)=\nobreak1$.
% In particular each such group is non-abelian.
\end{enumerate}
Moreover, for any such $G$ and $H$ we have:
\begin{enumerate}
\item[(c)] $G$ has Zassenhaus type (3,3) if and only if $\dim_{\BF_3}\gr_3^D(H)=2$.
\item[(d)] $H$ comes from a unique $\sigma$-isomorphism class of $\sigma$-$3$-groups generated by odd elements.
\item[(e)] The subgroup of $\sigma$-automorphisms $\Aut_\sigma(H)$ of $\Aut(H)$ has index~$3^2$.
\end{enumerate}
\end{Prop}

\begin{Proof}
The remarks above directly show the implication (a)$\Rightarrow$(b). For the converse note first that $\dim_{\BF_3}\gr_1^D(H)=2$ is equivalent to $d(H)=2$. From \eqref{F2Grad} we thus see that (b) is equivalent to $H$ being isomorphic to the quotient of $F_2/D_4(F_2)$ by a subspace of dimension $\le2$ of $\gr^D_3(F_2)$. Since $\gr^D_3(F_2)$ is totally odd, any such subspace is generated by two odd elements. As the formation of the odd part is exact in short exact sequences of $\sigma$-pro-$p$-groups by \cite[Prop.\,3.3]{Pink2025}, these elements can be lifted to odd elements of $F_2$. Letting $N\triangleleft F_2$ be the $\sigma$-invariant normal subgroup that is generated by all their conjugates, we thus deduce that $H\cong F_2/D_4(F_2)N \cong G/D_4(G)$ for the weak Schur $\sigma$-group $G:= F_2/N$, proving the implication (b)$\Rightarrow$(a).

Next, (c) follows at once from the definition of Zassenhaus type (3,3), and (d) is a special case of Proposition \ref{IsomVSSigmaIsom}.

To prove (e) we abbreviate $\bar F_2:=F_2/D_4(F_2)$ and set $A:=\Aut(\bar F_2)$ and $A_\sigma := \Aut_\sigma(\bar F_2)$. Also let $B$ denote the kernel of the homomorphism $A\onto\Aut(\gr^D_1(\bar F_2)) \cong \GL_2(\BF_3)$. As the restriction to $A_\sigma$ of this homomorphism is still surjective by \cite[Prop.\,6.2]{PinkRubio2025}, it follows that $A=BA_\sigma$. In particular this implies that $[A:A_\sigma] = [BA_\sigma:A_\sigma]=[B:B\cap A_\sigma]$. 

Pick odd generators $x$ and $y$ of $\bar F_2$. Then $B$ consists of all automorphisms $\phi$ of $\bar F_2$ with $\phi(x)\equiv x$ and $\phi(y)\equiv y$ modulo $D_2(\bar F_2)$. As any choice of elements $\phi(x)$, $\phi(y)\in\bar F_2$ satisfying these congruences results in an automorphism of~$\bar F_2$, it follows that $B$ has order $|D_2(\bar F_2)|^2$. On the other hand, $\phi$ is a $\sigma$-automorphism if and only if $\phi(x)$ and $\phi(y)$ are odd. By \cite[Prop.\,3.3]{PinkRubio2025} it follows that $B\cap A_\sigma$ has order $|D_2(\bar F_2)^-|^2$. But since $\gr_2^D(\bar F_2)$ is totally even of dimension~$1$ and $\gr_3^D(\bar F_2)$ is totally odd, we have $D_2(\bar F_2)^- = D_3(\bar F_2)$ and conclude that $[B:B\cap A_\sigma] = |\gr_2^D(\bar F_2)|^2 = 3^2$.

Now write $H$ as the quotient of $\bar F_2$ by a subspace $V\subset\gr_3^D(\bar F_2)$. 
Then every automorphism of $H$ lifts to an automorphism of $F_2$ and hence to an automorphism of $\bar F_2$; so we have a surjective homomorphism $\pi\colon\Stab_A(V)\onto \Aut(H)$. Likewise every $\sigma$-automorphism of~$H$ lifts to a $\sigma$-automorphism of $F_2$ by \cite[Prop.\,6.3]{PinkRubio2025}; so we have a surjective homomorphism $\pi_\sigma\colon\Stab_{A_\sigma}(V)\onto \Aut_\sigma(H)$. 
On the other hand, since $H$ surjects to $\bar F_2/D_3(\bar F_2)$, each automorphism $\phi\in\Ker(\pi)$ satisfies $\phi(x)\equiv x$ and $\phi(y)\equiv y$ modulo $D_3(\bar F_2)$. As $x$ and $y$ are odd, and $D_3(\bar F_2) = \gr_3^D(\bar F_2)$ is totally odd, from \cite[Prop.\,3.3]{PinkRubio2025} it follows that $\phi(x)$ and $\phi(y)$ are automatically odd as well. Thus $\phi$ lies in $A_\sigma$, and we deduce that $\Ker(\pi)=\Ker(\pi_\sigma)$. Combining this information we conclude that $[\Aut(H):\Aut_\sigma(H)] = [\Stab_A(V):\Stab_{A_\sigma}(V)]$. 

Now observe that $B$ acts trivially on $\gr_3^D(\bar F_n)$ by Proposition~\ref{AutFn}. In particular it stabilizes~$V$ and we have $B\cap\Stab_{A_\sigma}(V) = B\cap A_\sigma$. On the other hand the equality $A=BA_\sigma$ now implies that $\Stab_A(V) = B\cdot\Stab_{A_\sigma}(V)$. Together we deduce that 
$$\bigl[\Stab_A(V):\Stab_{A_\sigma}(V)\bigr]\ =\ 
%\bigl[B\cdot\Stab_{A_\sigma}(V):\Stab_{A_\sigma}(V)\bigr]\ =\ 
\bigl[B:B\cap\Stab_{A_\sigma}(V)\bigr]\ =\ 
\bigl[B:B\cap A_\sigma\bigr]\ =\ 3^2,$$
and are done.
\end{Proof}

%%%%%%%%%%%%%%%%%%%%%%%%
\medskip
Using the Small Groups Library built into the computer algebra system GAP, one can easily find all isomorphism classes of groups $H$ as in Proposition~\ref{HisGbyD4} (see Part~4 of our worksheet \cite[{\tt 1-Checking-powerfulness}]{Ahlqvist-Pink-I}).
The groups $G/D_4(G)$ for $G$ of Zassenhaus type $(3,3)$ have order $3^5=243$, and among all 67 isomorphism classes of groups of this order, these are the 13 groups obtained by the command {\tt SmallGroup(243,i)} for the indices
\UseTheoremCounterForNextEquation
\begin{equation}\label{GbyD4Type33Indices}
i\ \in\ \bigl\{2,\, 3,\, 4,\, 5,\, 6,\, 7,\, 8,\, 9,\, 13,\, 14,\, 15,\, 17,\, 18\bigr\}.
\end{equation}
The remaining groups are the 5 groups of order $3^6=729$ obtained by {\tt SmallGroup(729,i)} for the indices
\UseTheoremCounterForNextEquation
\begin{equation}\label{GbyD4Type3>3Indices}
i\ \in\ \bigl\{9,\ 10,\ 11,\ 12,\ 26\bigr\},
\end{equation}
and the group $F_2/D_4(F_2)$ of order $3^7=2187$ obtained by {\tt SmallGroup(2187,33)}.
In the following we call a weak Schur $\sigma$-group $G$ \emph{of type $[d,i]$} if $G/D_4(G)$ is isomorphic to the group {\tt SmallGroup(d,i)}. 
For $i$ in \eqref{GbyD4Type33Indices} we abbreviate the group {\tt SmallGroup(243,i)} by~$H_i$.

%%%%%%%%%%%%%%%%%%%%%%%%%%%%%%%%%%%%%%%%%%%%%%%%
%\newpage
\section{Analytic Schur $\sigma$-groups of type $(3,3)$ for $p=3$}
%\section{Powerful subgroups of small index}
\label{FrattPowerful}

For all weak Schur $\sigma$-groups $G$ of type $(3,3)$ with $p=3$, we now study which small index subgroups $E(G)$ are powerful. The smallest index subgroups relevant here are the first steps in the Zassenhaus filtration $D_k(F_2)$, yielding $E(G)=D_k(G)$ for $k\le4$. Observe that $D_1(G)=G$ is never powerful, because by Proposition \ref{HisGbyD4} we always have $\Agemo_1(G)\subset D_3(G)\subsetneqq D_2(G) = \Fratt(G)$ in this case.

\medskip
The GAP worksheet \cite[{\tt 1-Checking-powerfulness}]{Ahlqvist-Pink-I} finds that for 10 of the 13 cases in the list \eqref{GbyD4Type33Indices}, namely all those with the indices
\UseTheoremCounterForNextEquation
\begin{equation}\label{GbyD4Type33IndicesGood}
i\ \in\ \bigl\{2,\, 4,\, 5,\, 6,\, 7,\, 8,\, 14,\, 15,\, 17,\, 18\bigr\},
\end{equation}
the subgroup $D_2(G)$ is powerful with $d(D_2(G))\le3$ for every weak Schur $\sigma$-group $G$ of type $[243,i]$. In the remaining three cases $i\in\{3,9,13\}$ the group $D_2(G)$ is never powerful. The worksheet also gives partial results for the subgroups $D_3(G)$ and $D_4(G)$.

\medskip
The result for $D_2(G)$ has the following theoretical consequences:

\begin{Prop}\label{GoodTypePGLNullSet}
For each index $i$ in \eqref{GbyD4Type33IndicesGood}, we have:
\begin{enumerate}
\item[(a)] Any infinite strong Schur $\sigma$-group of type $[243,i]$ is isomorphic to an open subgroup of a form of $\PGL_2$ over~$\BQ_p$.
%The group $H$ from Bartholdi-Bush \cite{BartholdiBush2007} is an example for this.
\item[(b)] The set of $\sigma$-isomorphism classes of all infinite weak Schur $\sigma$-groups of type $[243,i]$ is a closed subset of~$\Sch$ of measure~$0$.
\end{enumerate}
\end{Prop}

\begin{Proof}
The computation shows that for any weak Schur $\sigma$-group of type $[243,i]$, the subgroup $D_2(G)$ is powerful with $d(D_2(G))\le3$. Thus (a) results from Proposition \ref{AnalStrongDim3PGL2}. Also, the set of $\sigma$-isomorphism classes of weak Schur $\sigma$-groups $G$ with $G/D_4(G)\cong H_i$ is open and closed in $\Sch$ by \cite[(10.1)]{PinkRubio2025}. Thus (b) follows from Proposition \ref{AnalDim3WeakSchurNullSet}
\end{Proof}

%%%%%%%%%%%%%%%%%%%%%%%%
\medskip
The  rest of this section is devoted to explaining the computations in more detail.

\medskip
To begin with, the computer algebra system GAP and especially its ANUPQ package have very efficient code to deal with finite $p$-groups, but none for $\sigma$-$p$-groups. While in principle it should be possible to adapt that code to $\sigma$-$p$-groups with the same level of efficiency, we shunned that effort, not being specialists. As our computation only looks at isomorphism classes of $\sigma$-$p$-groups, we use Proposition \ref{IsomVSSigmaIsom} instead. We therefore need a macro to decide whether a given finite $3$-group $G$ with $d(G)=2$ possesses an automorphism of order $2$ that acts as $-1$ on $G/\Fratt(G)$. The macro {\tt IsOurSigmaGroup} in Part~3 of our worksheet does that, using the inbuilt function {\tt AutomorphismGroupPGroup} that computes the automorphism group of a $p$-group. For large $G$ this is quite inefficient and possibly a limiting factor in our approach.

% NOTE FOR US: To compute the descendants of a $p$-group~$G$, the ANUPQ package uses the group $\Aut(G)$ and its orbits to isolate representatives for the isomorphism classes of descendants. For some groups of $p$-class 6 I compared the runtimes for computing all immediate descendants if $\Aut(G)$ has been computed before respectively not, and there does not seem to be a difference. So GAP does not seem to remember $\Aut(G)$. 
% (But maybe that can be done using Section 5.7-20 or 5.9-1 of the ANUPQ manual on ``interactive ANUPQ functions''?)
% Also, for the groups considered the computation of all immediate descendants used about 13 times as much time as the computation of $\Aut(G)$. So that would not make much difference. But restricted to computing only the immediate descendants with step size 3 took only about 4 times as long as computing $\Aut(G)$.]}

\medskip
The crucial point in our computation is that by Proposition \ref{EGisPowerful}, whether $E(G)$ is powerful or not depends only on the finite quotient $G/E_2(G)$. It therefore suffices to go through the groups $G/E_2(G)$ for all weak Schur $\sigma$-groups $G$ of type (3,3). When the algorithm determines that $E(G)$ is powerful, it explicitly constructs $G/E_2(G)$ up to isomorphism. When it determines that $E(G)$ is not powerful, it sometimes detects that on some smaller quotient $G/P_j(G)E_2(G)$ and then stops without having to generate $G/E_2(G)$.

Since $E$ is contained in $D_2(F_2)$, one can show that $E_2$ is contained in $D_4(F_2)$ (see Part~2 of our worksheet). Thus $E_2(G) \subset D_4(G)$, and so $G/D_4(G)\cong H_i$ is a quotient of $G/E_2(G)$. The problem thus separates completely according to the type $[243,i]$ of~$G$.

\medskip
The construction of $G/E_2(G)$ uses O'Brien's $p$-group generation algorithm  \cite{OBrien1990} implemented in the ANUPQ package. Based on the descending $p$-central series of~$G$, it constructs finite $p$-groups by induction on their $p$-class. Our algorithm thus inductively constructs all possibilities for $G/P_j(G)E_2(G)$ for $j\ge3$ as far as necessary.

\medskip
For $p=3$, the descending $p$-central series is connected to the Zassenhaus filtration by the relations
\UseTheoremCounterForNextEquation
\begin{equation}\label{DPDPDPInclusions}
D_2(G) = P_1(G) \supset D_3(G) \supset P_2(G) \supset D_4(G) \supset P_3(G).
\end{equation}
%These are direct consequences of the recursion formula \eqref{Zass} for $D_i(G)$ and that for $P_j(G)$. Indeed, for arbitrary $p$ the inclusions $D_i(G) \supset P_{i-1}(G)$ follow directly by induction over~$i$. The inclusion $P_2(G) \supset D_4(G)$ for $p=3$ results from the computation
%\begin{eqnarray*}
%D_4(G) &\stackrel{\eqref{Zass}}{=}& \Agemo_1(D_2(G))\cdot[D_1(G),D_3(G)]\cdot[D_2(G),D_2(G)] \\
%&\subset& \Agemo_1(D_2(G))\cdot[G,D_2(G)] \\
%&\subset& \Agemo_1(P_1(G))\cdot[G,P_1(G)] \\
%&=& P_2(G).
%\end{eqnarray*}
Moreover, all these inclusions are proper in the case $G=F_2$ (see Part~2 of our worksheet). 
The last inclusion in \eqref{DPDPDPInclusions} implies that $G/P_3(G)$ is an extension of $G/D_4(G)\cong H_i$. 

As initialization, for each $i$ we therefore create a list of all groups of the form $G/P_3(G)$ up to isomorphism for all weak Schur $\sigma$-groups $G$ of type $[243,i]$. This uses Proposition \ref{HKGEDProp} for $P_3(F_n)\subset D_4(F_n)$ in place of $E\subset D$ and is done in Parts 5--6 of our worksheet. 
%NOTE FOR US: Here $G/P_3(G)$ is not necessarily a descendant of $G/D_4(G)$. But it is a descendant of $G/P_2(G)$; so we can apply the ANUPQ function {\tt PqDescendants} to first compute all descendants of $G/P_2(G)$ of $p$-class $\le3$ and then identify the correct ones. This is reasonably fast. It might be faster using Proposition \ref{HKGDstarEDProp} or~\ref{HKGDstarPDProp}, but we didn't care that much.
From this list we can easily compute a list of possibilities for $G/P_3(G)E_2(G)$.
For later use, in Part 7 we also compute all possibilities for $G/P_4(G)$ and $G/P_5(G)$.

\medskip
The passage from $H=G/P_j(G)E_2(G)$ to $G/P_{j+1}(G)E_2(G)$ uses Proposition \ref{HKGDstarEDProp} for the subgroups $P_{j+1}(F_n)E_2\subset P_j(F_n)E_2$ in place of $E\subset D$. The assumptions of that proposition are satisfied, because we have $(P_j(F_n)E_2)^* = P_j(F_n)^*E_2^* = P_{j+1}(F_n)E_2^* \subset P_{j+1}(F_n)E_2$ by Proposition \ref{MstarProperties}~(a). Also, the assumption $P_j(H)=1$ implies $P_{j+1}(H^*)=1$. The relevant groups $G/P_{j+1}(G)E_2(G)$ up to isomorphism are therefore precisely the descendants $K$ of~$H$ which have $p$-class $\le j+1$ and satisfy
\begin{enumerate} 
\item[(a$^*$)] We have $E_2(K)=1$.
\item[(b$^*$)] The group $K$ possesses the structure of a $\sigma$-$p$-group that is generated by odd elements.
\item[(c$^*$)] We have $\mu(H) = \dim_{\BF_p}\Ker(K\!\onto\!H) + 2 + \dim_{\BF_p} E_2(H^*)$.
\end{enumerate}
To generate all these we first compute $\mu(H)-\dim_{\BF_p} E_2(H^*)-2$, apply the ANUPQ function {\tt PqDescendants} to $H$ with that step size according to (c$^*$), and then test for conditions (a$^*$) and (b$^*$). The code for this computation is contained in Part~8 of the worksheet. That part also contains macros for testing recursively, given an integer $j\ge3$ and a finite $p$-group $H$ of the form $H=G/P_j(G)$, whether the subgroup $E(G)$ is powerful for all such~$G$, respectively for no such~$G$, or whether both cases occur.
%using the condition from Proposition \ref{EGisPowerful} 
When $E(G)$ is powerful for all such~$G$, the code also computes the maximal rank of $E(G)$ in all those cases.

\medskip
In Part~9 this code is applied to the case $E=D_2$. As stated above, it finds that for 10 of the 13 cases in the list \eqref{GbyD4Type33Indices}, namely all those with the indices in \eqref{GbyD4Type33IndicesGood}, the subgroup $D_2(G)$ is powerful for every weak Schur $\sigma$-group $G$ of type $[243,i]$. Moreover in those cases it finds that the minimal number of generators of $D_2(G)$ is always~$\le3$. In the remaining 3 cases the group $D_2(G)$ is never powerful. 

\medskip
The rest of the worksheet computes partial results for $D_3(G)$ and $D_4(G)$ in the cases when $D_2(G)$ is not powerful. In Part~10 we find that for type $[243,13]$ the group $D_3(G)$ is never powerful, while in each of the types $[243,3]$ and $[243,9]$ both cases occur. 
In the cases that neither of $D_2(G)$ and $D_3(G)$ is powerful, we attempt to repeat the computation for $D_4(G)$ in Part~11. For the subcases of type $[243,9]$ we find that again both cases occur here. The types $[243,3]$ and $[243,13]$ were not completed due to memory overload. 

%Note for us: To find the $p$-class I compute the whole $p$-central series using the command {\tt PCentralSeries}. Can one use instead the function {\tt PqPClass} from Section 5.5-4 of the ANUPQ manual on ``interactive ANUPQ functions''? And would that be faster? (It's only a minor issue; most of the computation is spent in other parts.)

%%%%%%%%%%%%%%%%%%%%%%%%%%%%%%%%%%%%%%%%%%%%%%%%
%\newpage
\section{$p$-Tower groups of imaginary quadratic fields}
\label{pTG}

Now consider an imaginary quadratic field $K$ and let $\widetilde K$ denote its maximal unramified pro-$p$-extension for an odd prime~$p$. Then $\widetilde K$ is Galois over~$\BQ_p$, and its Galois group is the semidirect product of the $p$-tower group $G_K := \Gal(\widetilde K/K)$ with a group of order~$2$ generated by complex conjugation. Thus $G_K$ is a $\sigma$-pro-$p$-group. By Koch-Venkov \cite[\S1]{KochVenkov1974} and \cite[Prop.\,11.2]{PinkRubio2025} this is a strong Schur $\sigma$-group in the sense of Section~\ref{SigmaAndSchur}. 
In \cite[\S12]{PinkRubio2025} we explain that the probability measure $\mu_\infty$ on $\Sch$ is meant as a heuristic to describe statistical properties of $p$-tower groups of imaginary quadratic fields.

%%%%%%%%%%%%%%%%%%%%%%%%
\medskip
A particularly interesting question is when and how often $G_K$ is finite. Though much is known about that, the problem is not yet settled completely. 
According to a conjecture of McLeman \cite[Conj.\,2.9]{McLeman2008}, the $p$-tower group of an imaginary quadratic field should be finite if and only if it has $d(G_K)\le1$ or Zassenhaus type (3,3).
While the ``only if'' part of that conjecture has been refuted by the first author and Carlson \cite{AhlqvistCarlson2025}, the ``if'' part remains open.
The paper \cite{Pink2025} by the second author yields some credence for this in the case $p>3$.
In the case $p=3$ of the present paper, Proposition \ref{GoodTypePGLNullSet} has the same consequences for $p$-tower groups of type $[243,i]$ for any index $i$ in the set \eqref{GbyD4Type33IndicesGood}.

\medskip
Indeed, if the heuristic of \cite{PinkRubio2025} is correct, Proposition \ref{GoodTypePGLNullSet} (b) implies that McLeman's statement holds on average, namely that the proportion of finite groups $G_K$ among those of type $[243,i]$ tends to~$1$. Of course, this does not rule out a sparse infinite set of exceptions.

\medskip
But recall that the Fontaine-Mazur Conjecture \cite[Conj.\,5b]{FontaineMazur1995} predicts that $G_K$ does not have an infinite $p$-adic analytic quotient. In particular this implies that $G_K$ itself cannot be an infinite $p$-adic analytic group (compare Hajir \cite[Conj.\,2]{Hajir1997}). But if $G_K$ is of type $[243,i]$, Proposition \ref{GoodTypePGLNullSet} (a) does prove that $G_K$ is a $p$-adic analytic group. The Fontaine-Mazur conjecture therefore directly implies that $G_K$ is finite, yielding the ``if'' part of McLeman's conjecture in full for groups of type $[243,i]$.

%%%%%%%%%%%%%%%%%%%%%%%%%%%%%%%%%%%%%%%%%%%%%%%%
%\newpage
\section{Experimental results}
\label{Exp}

In this section we again assume that $p=3$. In this case the Cohen--Lenstra heuristic is supported by the seminal work of Davenport and Heilbronn \cite{DavenportHeilbronn1971}, concerning the average size of the $3$-torsion of class groups of quadratic fields.

%%%%%%%%%%%%%%%%%%%%%%%%
In Section \ref{33ByD4} we have seen that, for weak Schur $\sigma$-groups $G$ with $d(G)=2$, there are precisely 19 possibilities for the quotient $H:=G/D_4(G)$ up to isomorphism. The heuristic from \cite[\S12]{PinkRubio2025} predicts that each of these occurs with a specific asymptotic frequency, as $G$ runs through the $p$-tower groups $G_K$ for all imaginary quadratic fields $K$ with some bound on the discriminant.
By \eqref{MuUDG} the expected probability for a given $H$ is 
\UseTheoremCounterForNextEquation
\begin{equation}\label{muH1}
\mu(H)\ :=\ \frac{C_\infty}{C_{2-m}}\cdot\frac{1}{|\!\Aut_\sigma(H)|}\,,
\end{equation}
where $m=r_{D_4(F_2)}(H)$ is the minimal number of relations defining $H$ as a quotient of~$F_2/D_4(F_2)$. Thus 
\UseTheoremCounterForNextEquation
\begin{equation}\label{muHm}
m\ =\ \scriptstyle\left\{\displaystyle\begin{array}{ll}
2 & \hbox{for the 13 groups of type $[243,i]$ from \eqref{GbyD4Type33Indices},}\\[3pt]
1 & \hbox{for the 5 groups of type $[729,i]$ from \eqref{GbyD4Type3>3Indices},}\\[3pt]
0 & \hbox{for the group $F_2/D_4(F_2)$ of type $[2187,33]$.}
\end{array}\right.
\end{equation}
Using Proposition \ref{HisGbyD4} (e) this simplifies to
\UseTheoremCounterForNextEquation
\begin{equation}\label{muH2}
\mu(H)\ =\ 
\frac{C_\infty}{C_{2-m}}\cdot\frac{3^2}{|\!\Aut(H)|}\,,
\end{equation}
which can easily be computed by GAP.

On the other hand, for $p=3$ the Cohen--Lenstra heuristic for the maximal abelian quotient $G_{K,\ab}$ is known to converge quite slowly \cite{FungWilliams1990}, \cite{BelabasBhargavaPomerance2010}. In particular, the actual frequency of the event $d(G_K)=2$ converges only slowly to the expected probability from \eqref{MuSchn}:
\UseTheoremCounterForNextEquation
\begin{equation}\label{MuSch2}
\mu_\infty(\Sch_2)\ =\ \frac{C_\infty}{C_2^2} \cdot 3^{-4}\ \approx\ 0.01969.
\end{equation}
For this reason we focus on the conditional probability that $G/D_4(G)\cong H$ under the assumption $d(G)=2$, which 
%by combining \eqref{muH2} and \eqref{MuSch2} 
comes out as
\UseTheoremCounterForNextEquation
\begin{equation}\label{MuHCond}
\mu_\cond(H)\ :=\ 
\frac{\mu(H)}{\mu_\infty(\Sch_2)}\ =\ 
\frac{C_2^2}{C_{2-m}}\cdot\frac{3^6}{|\!\Aut(H)|}\,.
\end{equation}

%%%%%%%%%%%%%%%%%%%%%%%%
We have tested these values computationally as follows.
We begin with a list of all 461925 imaginary quadratic fields $K$ with $|d_K| < 10^8$ for which the class group has $p$-rank~2, obtained from the  $L$-functions and Modular Forms Database \cite{LMFDB}.
For each such field we compute an explicit presentation of $G_K/D_4(G_K)$ using the program \cite[{\tt 2-Massey}]{Ahlqvist-Pink-I}, as explained in Section~\ref{Present}. 
% In the huge majority of cases this computation finishes in less than a second, but it can take much longer or end in a stack overflow when the class group of $L_x$ is large. This happens in particular when $L_x$ has an unusually large $p'$-class group, which heuristically should not affect $G_K$ at all. We therefore impose a time limit of five minutes on each case and put the remaining ones aside. At present ... this concerns ... of the ... relevant discriminants, or approximately 0.06\%.
% We believe that this omission does not affect the significance of the observed frequencies, because these remained very robust as the number of exceptions was successively reduced by code improvements and the use of more powerful computers.

\medskip
Since the computation seems fraught with opportunities for errors such as sign mistakes, we have checked it for consistency by running most of it again with a different choice of the basis $x_1,x_2$ of $H^1(X, \BZ/p\BZ)$, which gave a different presentation of $G_K/D_4(G_K)$, but always the same isomorphism class. On the other hand, our results match perfectly with the index-$p$-abelianization-data (IPAD) from Boston--Bush--Hajir shown in \cite[Table 1]{BostonBushHajir2017}, as we will explain below.

%%%%%%%%%%%%%%%%%%%%%%%%
\medskip

In this way we have determined the isomorphism class of $G_K/D_4(G_K)$ for a total of $N:=461925$ discriminants $d_K$ with $d(G_K)=2$. Let $n(H)$ denote the number of occurrences of the isomorphism class of $H$ and $\mu_\obs(H) := n(H)/N$ the observed relative frequency. Then our results are collected in the following table, rounded to 5 decimal digits:
\UseTheoremCounterForNextEquation
\begin{equation}\label{Experi}
\begin{tabular}{|c|c|c|r|c|c|}
\hline
{\huge\mathstrut}
\kern0pt$H$      & $H_{\text{ab}}$  & {$\mu_\cond(H)$} & $n(H)$ & {$\mu_\obs(H)$} & {$\mu_\obs(H)/\mu_{\rm cond}(H)$}\\[5pt] 
\hline\hline
{\Large\mathstrut}%
$[ 243, 2 ]$	 & $[ 9, 9 ]$	  & $0.00732$	  & $3184$	  & $0.00689$	  & $0.94175$ \\    [5pt]
$[ 243, 3 ]$	 & $[ 3, 3 ]$	  & $0.04392$	  & $19298$	  & $0.04178$	  & $0.95131$ \\ 
$[ 243, 4 ]$	 & $[ 3, 3 ]$	  & $0.08783$	  & $40968$	  & $0.08869$	  & $1.00977$ \\ 
$[ 243, 5 ]$	 & $[ 3, 3 ]$	  & $0.17566$	  & $83353$	  & $0.18045$	  & $1.02724$ \\ 
$[ 243, 6 ]$	 & $[ 3, 3 ]$	  & $0.08783$	  & $40125$	  & $0.08686$	  & $0.98899$ \\ 
$[ 243, 7 ]$	 & $[ 3, 3 ]$	  & $0.08783$	  & $41398$	  & $0.08962$	  & $1.02037$ \\ 
$[ 243, 8 ]$	 & $[ 3, 3 ]$	  & $0.08783$	  & $40807$	  & $0.08834$	  & $1.00580$ \\ 
$[ 243, 9 ]$	 & $[ 3, 3 ]$	  & $0.02196$	  & $10426$	  & $0.02257$	  & $1.02791$ \\     [5pt]
$[ 243, 13 ]$	 & $[ 3, 9 ]$	  & $0.02928$	  & $13288$	  & $0.02877$	  & $0.98256$ \\   
$[ 243, 14 ]$	 & $[ 3, 9 ]$	  & $0.02928$	  & $13705$	  & $0.02967$	  & $1.01340$ \\   
$[ 243, 15 ]$	 & $[ 3, 9 ]$	  & $0.02928$	  & $13474$	  & $0.02917$	  & $0.99632$ \\   
$[ 243, 17 ]$	 & $[ 3, 9 ]$	  & $0.08783$	  & $39425$	  & $0.08535$	  & $0.97174$ \\   
$[ 243, 18 ]$	 & $[ 3, 9 ]$	  & $0.17566$	  & $81494$	  & $0.17642$	  & $1.00433$ \\     [3pt]
\hline
{\Large\mathstrut}%
$[ 729, 9 ]$	  & $[ 3, 9 ]$	  & $0.00488$	  & $1979$	  & $0.00428$	  & $0.87801$ \\ 
$[ 729, 10 ]$	  & $[ 3, 9 ]$	  & $0.01464$	  & $6555$	  & $0.01419$	  & $0.96940$ \\    
$[ 729, 11 ]$	  & $[ 3, 9 ]$	  & $0.01464$	  & $6172$	  & $0.01336$	  & $0.91276$ \\    
$[ 729, 12 ]$	  & $[ 3, 9 ]$	  & $0.00976$	  & $4299$	  & $0.00931$	  & $0.95365$ \\  [5pt]
$[ 729, 26 ]$	  & $[ 9, 9 ]$	  & $0.00488$	  & $1929$	  & $0.00418$	  & $0.85582$ \\  [3pt]
\hline
{\Large\mathstrut}%
$[ 2187, 33 ]$	  & $[9,9]$	  & $0.00015$	  & $46$	  & $0.00010$	  & $0.65307$ \\  [3pt]
\hline
\end{tabular}
\end{equation}
Since the values of $\mu_\cond(H)$ vary over several orders of magnitude, we view the agreements in the last column as a moderate confirmation of the heuristic. 

\medskip
The worst value in the last column is that for $H = [2187,33]$. That is the case where all triple Massey products vanish, which by \cite[Thm.\,4.22]{AhlqvistCarlson2025} happens if and only if all unramified cyclic extensions of degree $p$ of $K$ have class group of $p$-rank at least 4. Thus that case is related to large $p$-class groups, and so the bad fit might be explained by the slow convergence of the Cohen--Lenstra heuristic for $p=3$ that was already mentioned above.
% We did use \cite[Thm.\,4.22]{AhlqvistCarlson2025} to check that none of the omitted cases is of type $[2187,33]$, so the discrepancy is real.
%not due to the omission of cases due to overflow.
%approximately 0.2\% missing data. 

%%%%%%%%%%%%%%%%%%%%%%%%
\medskip
After finishing these computations, we compared our results with the index-$p$-abelian\-iza\-tion-data (IPAD) from Boston--Bush--Hajir \cite[\S4]{BostonBushHajir2017}. Recomputing the IPAD for all 461925 discriminants we found that, in all those cases, the IPAD determines the isomorphism class of $G_K/D_4(G_K)$; see \cite[{\tt 5-Data-analysis/IPADs.txt}]{Ahlqvist-Pink-I}. Thus for the 14 most frequent IPADs shown in \cite[Table~1]{BostonBushHajir2017}, that quotient has type $[243, n]$ with $n= 5$, $18$, $7$, $4$, $8$, $6$, $18$, $3$, $17$, $9$, $14$, $13$, $17$,~$15$, respectively. Moreover, precisely 4 of the 19 types from Table \eqref{Experi} correspond to a unique IPAD,  namely:
       \begin{equation}
              \begin{tabular}{|c|c|c|c|}
              \hline
              {\Large\mathstrut}
              \kern0pt$H$      & IPAD   & $n(H)$ & $n(\text{IPAD})$\\[2pt] 
              \hline\hline
              {\Large\mathstrut}%
              $[ 243, 4 ]$	 & $[3, 3]; [3, 3, 3]^3, [9, 3]$		  & $40968$ & $40968$ \\ 
              $[ 243, 5 ]$	 & $[3, 3]; [3, 3, 3], [9, 3]^3$		  & $83353$ & $83353$ \\
              $[ 243, 7 ]$	 & $[3, 3]; [3, 3, 3]^2, [9, 3]^2$		  & $41398$ & $41398$ \\
              $[ 243, 9 ]$	 & $[3, 3]; [9, 3]^4$		         & $10426$ & $10426$  \\    
              \hline
              \end{tabular}
              \end{equation}
Since we used a very different approach than \cite{BostonBushHajir2017}, but arrived at the same numbers and matching lists of discriminants, this gives very strong evidence that both their and our computations are correct.

\medskip
The outcome suggests that the IPAD might always determine the isomorphism class of $G_K/D_4(G_K)$. In other words, it suggests that for any strong Schur $\sigma$-group $G$ with ${d(G)=2}$, the isomorphism class of its abelianization $G_\ab$ together with the isomorphism classes of $H_\ab$ for all subgroups $H\subset G$ index~$3$ might determine the isomorphism class of $G/D_4(G)$. But we have not 
%been able to decide 
settled this matter in general. 

%%%%%%%%%%%%%%%%%%%%%%%%
\medskip
However, we have found that the IPAD does \emph{not} determine the isomorphism class of $G/P_3(G)$ in general. For example, there exist strong Schur $\sigma$-groups $G$ with equal IPADs $[ 3, 9 ]; [ 3, 3, 9 ], [ 3, 27 ]^3$ and $G/D_4(G)$ of type $[243, 14]$, but with distinct isomorphism classes of $G/P_3(G)$, namely of types $[729, 17]$ and $[729, 20]$ in the Small Groups Library from GAP~\cite{GAP4}.

\medskip
On the other hand, we do have $G/D_4(G)\cong G/P_3(G)$ whenever the abelianization $G_\ab$ has elementary divisors $[3,3]$, that is, whenever $G/D_4(G)$ is of type $[243,n]$ with $3\le n\le 9$. For each such type we have thus determined the number of occurences of that isomorphism class of $G_K/P_3(G_K)$ among all imaginary quadratic fields $K$ with $|d_K|<10^8$.

%%%%%%%%%%%%%%%%%%%%%%%%%%%%%%%%%%%%%%%%%%%%%%%%
%\newpage

%%%%%%%%%%%%%%%%%%%%%%%%%%%%%%%%%%%%%%%%%%%%%%%%


\begin{thebibliography}{99}
% \renewcommand{\baselinestretch}{0.5}\normalsize
\setlength{\parskip}{0pt}
\addcontentsline{toc}{section}{References}
\footnotesize

\bibitem{AhlqvistCarlson2023}
Ahlqvist, E.; Carlson, M.:
{\rm The \'etale cohomology ring of the ring of integers of a number field,}
{\it Res. number theory} (2023) 9:58.
\url{ https://doi.org/10.1007/s40993-023-00467-5}

\bibitem{AhlqvistCarlsonOpen}
Ahlqvist, E.; Carlson, M.:
{\rm The \'etale cohomology ring of a punctured arithmetic curve,}
arXiv:2110.01597 [math.NT], 2021.
\url{ https://arxiv.org/abs/2110.01597}

\bibitem{AhlqvistCarlson2025}
Ahlqvist, E.; Carlson, M.:
{\rm Massey products in the \'etale cohomology of number fields,}
{\it J. reine angew. Math.} {\bf 823} (2025), 61--112.
\url{ https://doi.org/10.1515/crelle-2025-0006}

\bibitem{Ahlqvist-Pink-I}
Ahlqvist, E.; Pink, R.:
{\it Worksheets associated with this article,} 
\url{https://github.com/ericahlqvist/Ahlqvist-Pink-I}

\bibitem{ANUPQ}
Gamble, G; Nickel, W; O'Brien, E.:
{\rm ANUPQ, ANU p-Quotient,}
Version 3.3.3, GAP package, Nov.\ 2025.
\url{ https://gap-packages.github.io/anupq/}.

\bibitem{ArtinVerdierSeminar1964}
Artin, E.; Verdier, J.-L.:
{\it Seminar on \'etale cohomology of number fields,}
Lecture Notes, Woods Hole Summer Institute on Algebraic Geometry, 1964.

\bibitem{BartholdiBush2007}
Bartholdi, L.; Bush, M.:
{\rm Maximal unramified 3-extensions of imaginary quadratic fields and $\SL_2(\BZ_3)$,}
{\it J. Number Theory} {\bf 124} (2007), no. 1, 159--166.

\bibitem{BelabasBhargavaPomerance2010}
Belabas, K.; Bhargava, M.; Pomerance, C.:
{\rm Error estimates for the Davenport--Heilbronn theorems,}
{\it Duke Math. J.} {\bf 153} (2010), no. 1, 173--210.

\bibitem{Berkovich2008}
Berkovich, Y.:
{\it Groups of prime power order, Vol. 1,}
De Gruyter Exp. Math. {\bf 46}, Berlin 2008.

\bibitem{BostonBushHajir2017}
Boston, N.; Bush, M. R.; Hajir, F.:
{\rm Heuristics for $p$-class towers of imaginary quadratic fields,}
{\it Math. Ann.} {\bf 368} (2017) 633--669.

\bibitem{BostonNover2006}
Boston, N.; Nover, H.:
{\it Computing pro-$p$ Galois groups,}
Algorithmic number theory,
Lecture Notes in Computer Science 4076, Springer 2006, 1--10.

\bibitem{CohenLenstra1984}
Cohen, H.; Lenstra, H. W.:
{\it Heuristics on class groups of number fields,}
Number Theory, Noordwijkerhout 1983.
Lecture Notes in Math. 1068, Springer 1984, 33--62.
\url{ https://doi.org/10.1007/BFb0099440}

\bibitem{DavenportHeilbronn1971}
Davenport, H.; Heilbronn, H.:
{\rm On the density of discriminants of cubic fields. II,}
{\it Proc. Roy. Soc. Lond. A} {\bf 322} (1971), no. 1551, 405--420. 

\bibitem{DdSMS2003}
Dixon, J. D.; du Sautoy, M. P. F.; Mann, A.; Segal, D.:
{\it Analytic pro-$p$ groups,}
Cambridge University Press 1999.

\bibitem{DwyerMassey}
Dwyer, W. G.:
{\rm Homology, Massey products and maps between groups,}
{\it J. Pure Appl. Algebra} {\bf 6} (1975), 177--190.

\bibitem{FontaineMazur1995}
Fontaine, J.-M.; Mazur, B.:
{\it Geometric Galois representations,}
Elliptic curves, modular forms, \& Fermat's last theorem
(Hong Kong, 1993),
Ser. Number Theory, I, Int. Press, Cambridge, MA, 1995, 41--78.

\bibitem{FungWilliams1990}
Fung, G. W.; Williams, H. C.:
{\rm On the computation of a table of complex cubic fields with discriminant $D > -10^6$,}
{\it Math. Comp.} {\bf 55} (1990), no. 191, 313--325.

\bibitem{GAP4}
The GAP~Group, 
{\it GAP -- Groups, Algorithms, and Programming,}
Version 4.15.1,
2025,
\url{ https://www.gap-system.org}.

\bibitem{GolodShafarevich1964}
Golod, E. S.; Shafarevich, I. R.:
{\rm On the class field tower,}
{\it Izv. Akad. Nauk SSSR} {\bf 28} (1964), 261--272.

\bibitem{Hajir1997}
Hajir, F.:
{\rm On the Growth of $p$-Class Groups in $p$-Class Field Towers,}
{\it J. of Algebra} {\bf 188} (1997), 256--271.

\bibitem{Koch2002}
Koch, H.:
{\it Galois theory of $p$-extensions,}
Springer 2002.

\bibitem{KochVenkov1974}
Koch, H.; Venkov, B. B.:
{\rm The $p$-tower of class fields for an imaginary quadratic field,}
(Russian)
Modules and representations,
{\it Zap. Nau\v{c}n. Sem. Leningrad. Otdel. Mat. Inst. Steklov. (LOMI)} 46 (1974), 5--13.
English Translation: {\it Journal of Soviet Mathematics} {\bf 9} (1978), 291--299.

\bibitem{Lazard1965}
Lazard, M.:
{\it Groupes analytiques $p$-adiques,}
Inst. Hautes \'Etudes Sci. Publ. Math. {\bf 26} (1965), 389--603.

\bibitem{LMFDB}
The LMFDB Collaboration:
{\it The $L$-functions and Modular Forms Database,}
\url{ https://www.lmfdb.org/NumberField/QuadraticImaginaryClassGroups}, 2025.

\bibitem{MazurEtale}
Mazur, B.:
{\rm Notes on \'etale cohomology of number fields,}
{\it Ann. Sci. École Norm. Sup.} {\bf 6} (1973), 521--552.

\bibitem{McLeman2008}
McLeman, C.:
{\rm $p$-tower groups over quadratic imaginary number fields,}
{\it Ann. Sci. Math. Qu\'ebec} {\bf 32} (2008), no. 2, 199--209.

\bibitem{McLeman2009}
McLeman, C.:
{\rm A Golod--Shafarevich equality and $p$-tower groups,}
{\it J. Number Theory} {\bf 129} (2009), no. 11, 2808--2819.

\bibitem{MilneEtale}
Milne, J. S.:
{\it \'Etale Cohomology,}
Princeton Mathematical Series, vol.~33,
Princeton University Press, 1980.

\bibitem{NSW2008}
Neukirch, J.; Schmidt, A.; Wingberg, K.:
{\it Cohomology of Number Fields,}
Springer 2008.

\bibitem{OBrien1990}
O'Brien, E. A.:
{\rm The $p$-group generation algorithm,}
{\it J. Symbolic Comput.} {\bf 9} (5-6) (1990), 677--698.

\bibitem{PARI}
The PARI~Group, 
{\it PARI/GP}, 
version 2.18.1, Univ. Bordeaux, 2025.
\url{ http://pari.math.u-bordeaux.fr/}

\bibitem{Pink2025}
Pink, R.:
{\rm Schur $\sigma$-groups of type $(3,3)$,}
arXiv:2505.05580, [math.NT] (2025), 
\url{ https://arxiv.org/abs/2505.05580}

\bibitem{PinkRubio2025}
Pink, R.; Rubio, L. \'A.:
{\rm A Cohen--Lenstra Heuristic for Schur $\sigma$-Groups,}
arXiv:2505.05569, [math.NT] (2025), 
\url{ http://arxiv.org/abs/2505.05569}

\bibitem{RoquetteCFT}
Roquette, P.:
{\it On class field towers,}
Algebraic Number Theory, Proc.\ Instructional Conf., Brighton, 1965,
Academic Press, London, 1967, 231--249.

\bibitem{ScholzTaussky1934}
Scholz, A.; Taussky, O.:
{\rm Die Hauptideale der kubischen Klassenk\"orper imagin\"ar-quadratischer Zahlk\"orper,}
{\it J. Reine Angew. Math.} {\bf 171} (1934) 19--41.

\bibitem{Vogel-Massey}
Vogel, D.:
{\it Massey products in the Galois cohomology of number fields,}
PhD thesis,
Universität Heidelberg, February 10, 2004.
\url{ https://doi.org/10.11588/heidok.00004418.}

\end{thebibliography}
\end{document}